\documentclass[10pt,a4paper]{amsart}
\usepackage[arrow, matrix, curve]{xy}
\usepackage{amssymb}

\usepackage{color}
\usepackage{hyperref}

\usepackage{tikz}

\newtheorem{thm}{Theorem}
\newtheorem{lem}{Lemma}
\newtheorem{cor}{Corollary}

\newtheorem{pro}{Proposition}
\newtheorem{df}{Definition}
\newtheorem{conj}{Conjecture}

\newenvironment{pf}{{\it Proof:}\quad}{\hfill$QED$}
\newcommand{\R}{\mathbb R}
\newcommand{\N}{\mathbb N}
\newcommand{\C}{\mathbb C}
\newcommand{\Z}{\mathbb Z}
\newcommand{\Q}{\mathbb Q}
\newcommand{\PC}{\widehat{\mathcal{P}}}
\newcommand{\oz}{\overline{z}}
\newcommand{\oa}{\overline{a}}
\newcommand{\ob}{\overline{b}}
\newcommand{\oc}{\overline{c}}
\newcommand{\od}{\overline{d}}
\newcommand{\oee}{\overline{e}}
\newcommand{\of}{\overline{f}}

\begin{document}

\title{Fundamental classes of 3-manifold groups representations in SL(4,R)}
 \author{Thilo Kuessner}
   \address{School of Mathematics,
   KIAS, Hoegi-ro 85, Dongdaemun-gu,
   Seoul, 130-722, Republic of Korea}
   \email{kuessner@kias.re.kr}

\maketitle

\begin{abstract}We compute the fundamental class (in the extended Bloch group) for representations of fundamental groups of $3$-manifolds to $SL(4,\R)$ that factor over $SL(2,\C)$, in particular for those factoring over the isomorphism $PSL(2,\C)=SO(3,1)$. We also discuss consequences for the number of connected components of $SL(4,\R)$-character varieties, and we show that there are knots with arbitrarily many components of vanishing Chern-Simons invariant in their $SL(n,\C)$-character varieties.
\end{abstract}

\section{Introduction}
Given a finitely generated group $\Gamma$ and a connected Lie group $G$, it is a natural and fruitful question to understand the topology of its representation variety $Hom(\Gamma,G)$ and its character variety $X(\Gamma,G)=Hom(\Gamma,G)//G$, in particular to distinguish the connected components of these varieties.\footnotemark\footnotetext{For a manifold with boundary, we will consider the character variety of characters of boundary-unipotent representations $X_{bup}(\Gamma,G)$ only, which frequently has more connected components than $X(\Gamma,G)$, see e.g.\ \cite{mp}.} 

This question has been studied mostly for $\Gamma=\pi_1\Sigma_g$, the fundamental group of a closed, orientable surface of genus $g\ge 2$.
For compact Lie groups $G$ it is known since the work of Atiyah and Bott that the components of the character variety $X(\pi_1\Sigma_g,G)$ can be distinguished by the values of characteristic classes (of the flat bundles associated to the respective representations). This is however in general not true for noncompact Lie groups. It is known from work of Hitchin that for $n\ge 3$ the character variety $X(\pi_1\Sigma_g,PSL(n,\R))$ has $3$ or $6$ components, according to whether $n$ is odd or even, and not all of these components can be distinguished by characteristic classes. On the other hand, $X(\pi_1\Sigma_g,PSL(2,\R))$ has $4g-3$ components distinguished by the possible values of the Euler class, according to Goldman. 

Not much is known for the case that $\Gamma=\pi_1M$ is the fundamental group of a (closed or cusped) hyperbolic $3$-manifold. In this case, an invariant which can distinguish different components of $X_{bup}(\Gamma,SL(n,\C))$ is the Cheeger-Chern-Simons invariant of the flat bundle associated to a representation. This is a complex-valued invariant and it has become common to denote its imaginary part by the "volume" and the negative of its real part by the "Chern-Simons invariant" of the representation. These names stem from the case of hyperbolic metrics (more precisely, the lifts of their monodromy representation to $SL(2,\C)$), where these invariants yield indeed the volume and the Chern-Simons invariant of the hyperbolic metric.

It was proved in \cite{gtz} that the Cheeger-Chern-Simons invariant can be computed from the "fundamental class of $\rho$" 
in the so-called extended Bloch group $\widehat{\mathcal B}(\C)$, namely that it is the result of applying 
the Rogers' extended dilogarithm to that fundamental class. The group $\widehat{\mathcal B}(\C)$
is isomorphic to $H_3(SL(2,\C);\Z)$ and to a direct summand of $H_3(SL(n,\C);\Z)$ for $n\ge 3$. One of the main results 
in \cite{gtz} is to give a formula which computes that fundamental class 
from the naturally defined class $(B\rho)_*\left[M\right]\in H_3(SL(n;\C);\Z)$, which is the image of the fundamental class $\left[M\right]\in H_3(M;\Z)$ under the classifying map $B\rho\colon M\to BSL(n,\C)^\delta$. In the cusped case one has to consider $(B\rho)_*\left[M,\partial M\right]\in H_3(SL(n,\C),N;\Z)$ for $N\subset SL(n,\C)$ the subgroups of upper triangular matrices with 1's on the diagonal and \cite{gtz} also computes the fundamental class in the extended Bloch group from this. Because of this close connection we will denote the "fundamental class of $\rho$" in the extended Bloch group by 
$$\rho_*\left[M,\partial M\right]\in \widehat{\mathcal B}(\C).$$

In this paper we are going to consider finite-volume hyperbolic $3$-manifolds $M$ and representations $\rho\colon \pi_1M\to SL(4,\C)$ of their fundamental groups which factor over a representation $SL(2,\C)\to SL(4,\C)$. 
The main part of the paper will be devoted to the computation in the case of the 2-fold covering $SL(2,\C)\to SO(3,1)$ because this is the only irreducible $SL(2,\C)$-representation for which the computation does not already follow easily from the results in \cite{gtz}. 

{\bf Behaviour of fundamental classes under the isomorphism $PSL(2,\C)\to SO(3,1)$.}
\begin{thm}\label{thm1}Let $M$ be a compact, orientable $3$-manifold.
Let $\tau$ be the isomorphism $\tau\colon PSL(2,\C)\to SO(3,1)$ and let $\rho\colon\pi_1M\to PSL(2,\C)$ any boundary-unipotent representation. 
Then $$(\tau\circ\rho)_*\left[M,\partial M\right]=2\rho_*\left[M,\partial M\right]+
2\overline{\rho_*\left[M,\partial M\right]}\in\widehat{\mathcal{B}}(\C)$$
if $\rho$ lifts to a boundary-unipotent representation $\pi_1M\to SL(2,\C)$ 
(in particular if $M$ is closed)
and\footnotemark\footnotetext{The notation $\widehat{\mathcal{B}}(\C)_{PSL}$ will be explained in the last paragraph of \hyperref[ccsgtz]{Section \ref*{ccsgtz}}.}  
$$(\tau\circ\rho)_*\left[M,\partial M\right]=2\rho_*\left[M,\partial M\right]+
2\overline{\rho_*\left[M,\partial M\right]}\in\widehat{\mathcal{B}}(\C)_{PSL}$$
otherwise.
\end{thm}

This will be proved in \hyperref[compu]{Section \ref*{compu}}. Rather than computing $(\tau\circ\rho)_*\left[M,\partial M\right]$
for $\tau\colon PSL(2,\C)\to SO(3,1)$ it turns out to be more convenient to consider 
the equivalent representation $\rho\otimes\overline{\rho}\sim\tau\circ\rho$. (The problem with considering $\tau\circ\rho$ 
would be 
that a simplex with
an in the sense of \hyperref[ptolemy1]{Definition \ref*{ptolemy1}} generic $PSL(2,\C)/N$-decoration need not have a 
generic $SO(3,1)/N$-decoration, so further subdivision of the 
triangulation would be necessary to obtain generic $SO(3,1)/N$-decorations. )

Direct applications of the formulas from \cite{gtz} to $\rho\otimes\overline{\rho}$ will a priori yield some apparently unmanageable formula for $$(\tau\circ\rho)_*\left[M,\partial M\right]=(\rho\otimes\overline{\rho})_*\left[M,\partial M\right]\in \widehat{\mathcal{B}}(\C).$$ 
However, in the end almost everything in this formula will cancel out to yield the simple formula 
in \hyperref[thm1]{Theorem \ref*{thm1}}. The reason behind this will be cancellations in the extended 
pre-Bloch group\footnotemark\footnotetext{Remarkably all of our computations will work in the extended pre-Bloch group 
$\PC(\C)$ and the fact that the fundamental class actually belongs to the subgroup $\widehat{\mathcal B}(\C)\subset\PC(\C)$ 
will not be needed.}, using the 5-term relation and some symmetries
in $\PC(\C)$ that we describe in \hyperref[symmetries]{Section \ref*{symmetries}}. 

Our original motivation for computing the fundamental class of these representations was to compute their Chern-Simons invariants. Namely, it follows from \hyperref[thm1]{Theorem \ref*{thm1}} that $Vol(\tau\circ\iota)=0$ and\footnotemark\footnotetext{If $M$ has cusps, then this is of course to be understood as an equality modulo $\pi^2$ because $CS(M)$ is only defined up to this ambiguity.}
$CS(\tau\circ\iota)=4CS(M)$. It turns out, however, that this equality can be proved by an easier (and well-known) argument, which we will give in \hyperref[invpol]{Lemma \ref*{invpol}}. 

{\bf Geometric $SL(n,\C)$-representations of $3$-manifold groups.}
Recall that for an oriented, hyperbolic $3$-manifold $M$, its fundamental group $\pi_1M$ is a discrete subgroup of $Isom^+(H^3)=PSL(2,\C)$. If $M$ is closed, then by \cite[Corollary 2.2]{cul} it lifts to a discrete subgroup $\Gamma\subset SL(2,\C)$. (If $M$ has cusps, then this lift is in general not boundary-unipotent, see \hyperref[nonlift]{Section \ref*{nonlift}} for a discussion of this case.) We call $\iota\colon \pi_1M\to SL(2,\C)$ (the lift of) the hyperbolic monodromy.

For each natural number $n$, the irreducible representation $\rho_n\colon SL(2,\C)\to SL(n,\C)$ corresponding to the unique $n$-dimensional $\C$-linear representation of the Lie algebra $\mathfrak{sl}(2,\C)$, can be composed with $\iota$ to yield a representation $\rho_n\circ\iota$ that is called a "geometric representation". Garoufalidis-D.Thurston-Zickert used their methods to give a short and elegant proof for the formula
$$(\rho_n\circ\iota)_*\left[M,\partial M\right]=\left(\begin{array}{c}n+1\\
3\end{array}\right)\iota_*\left[M,\partial M\right]\in \widehat{\mathcal B}(\C),$$
which implies the equalities
$$Vol(\rho_n\circ\iota)=\left(\begin{array}{c}n+1\\
3\end{array}\right)Vol(\iota),\hspace{0.5in}
CS(\rho_n\circ\iota)=\left(\begin{array}{c}n+1\\
3\end{array}\right)CS(\iota),$$
see \cite[Theorem 11.3]{gtz}. Actually also these equalities can be derived from properties of invariant polynomials, see \hyperref[invpol]{Lemma \ref*{invpol}}.

In \hyperref[tensconj]{Conjecture \ref*{tensconj}} we propose that the natural generalization of \hyperref[thm1]{Theorem \ref*{thm1}} should be an equality $$(\rho_n\otimes\overline{\rho_m})_*\left[M,\partial M\right]=m\left[M,\partial M\right]+n\overline{\left[M,\partial M\right]}$$ for all $n,m\ge 2$.
  
{\bf 4-dimensional representations of 3-manifold groups.} The classification of representations of the Lorentz group implies that with the exception of $\rho_2\otimes\overline{\rho}_2$ all $4$-dimensional representations of $SL(2,\C)$ are, up to conjugacy in $GL(4,\C)$,  
obtained as direct sums of the geometric representations and their complex conjugates. The exceptional case $\rho_2\otimes\overline{\rho}_2$ is equivalent to the $2$-fold covering $SL(2,\C)\to SO(3,1)$.

The Cheeger-Chern-Simons invariant is additive under direct sum and takes the complex conjugate upon complex conjugation of the representation, see \hyperref[properties]{Section \ref*{properties}}. From these principles and the computation of volumes and Chern-Simons invariants for the $\rho_n$'s and $\rho_2\otimes\overline{\rho_2}$ we obtain the following table of volumes and Chern-Simons invariants for the representations $\rho\circ\iota\colon\pi_1M\to GL(4,\C)$, where $\iota\colon\pi_1M\to SL(2,\C)$ is a lift of the monodromy of a hyperbolic structure and $\rho$ runs over all representations $\rho\colon SL(2,\C)\to GL(4,\C)$.
\begin{center}\begin{tabular}{|l|l|l|}
\hline
representation $\rho$& Volume of $\rho\circ\iota$ &Chern-Simons invariant of $\rho\circ\iota$\\
\hline
$\rho_4$&10 Vol(M)&10 CS(M)\\
$\overline{\rho_4}$&-10 Vol(M)&10 CS(M)\\
$\rho_2\otimes\overline{\rho}_2$&0&4 CS(M)\\
$\rho_3\oplus 1$&4 Vol(M)&4 CS(M)\\
$\overline{\rho}_3\oplus 1$&-4 Vol(M)&4 CS(M)\\
$\rho_2\oplus\rho_2$&2 Vol(M)&2 CS(M)\\
$\rho_2\oplus\overline{\rho}_2$& 0& 2 CS(M)\\
$\rho_2\oplus 1\oplus 1$&Vol (M)&CS (M)\\
$\overline{\rho}_2\oplus 1\oplus 1$&-Vol (M)& CS(M)\\
$1^{\oplus4}$&0&0\\
\hline
\end{tabular}\end{center}
For a hyperbolic $3$-manifold with nonvanishing Chern-Simons invariant $CS(M)\not=0$ this implies that we can distinguish 10 components of the $SL(4,\C)$-character variety by volume and Chern-Simons invariant.\footnotemark\footnotetext{One obtains an analogous table of volumes and Chern-Simons invariants for $\rho\circ\kappa$ whenever $\kappa\colon\pi_1M\to SL(2,\C)$ is any given representation. So, when the $SL(2,\C)$-character variety already has several components, then in most cases one will find accordingly more components in the $SL(4,\C)$-character variety.} 

We remark that in the above list the three representations of vanishing volume can be conjugated into $SL(4,\R)$. Indeed, $\rho_2\otimes\overline{\rho}_2$ is equivalent to the well-known 2-fold covering map $SL(2,\C)\to SO(3,1)\subset SL(4,\R)$, while $\rho_2\oplus \overline{\rho_2}$ is equivalent to the embedding $SL(2,\C)\to SL(4,\R)$ coming from $(a_1+a_2i)\to\left(\begin{array}{cc}a_1&a_2\\
-a_2&a_1\end{array}\right)$, see \hyperref[4r]{Section \ref*{4r}}.
So we obtain the following corollary.
\begin{cor}\label{3compo}For finite-volume hyperbolic, orientable $3$-manifolds with nonvanishing Chern-Simons invariant there are at least 3 connected components in their $SL(4,\R)$-character variety.\end{cor}

{\bf More components.} The ptolemy module \cite{pto} in SnapPy \cite{snappy} computes the $SL(2,\C)$-representations for $3$-manifolds built from up to $9$ ideal tetrahedra and the $SL(3,\C)$-representations for many $3$-manifolds built from up to $4$ tetrahedra\footnotemark\footnotetext{Another approach to the computation of $SL(3,\C)$-representations for $3$-manifold groups uses the Fock-Goncharov coordinates, it is explained in \cite{fkr} and at the time of writing it computes them for a similar range than the ptolemy module.}, but at the time of writing can compute $SL(4,\C)$-representations only for $3$-manifolds composed of two ideal tetrahedra, that is for the figure eight knot complement and its sister. (For them it detects only irreducible representations because the reducible ones would need more than two simplices to allow a generic decoration.) It turns out that for the figure eight knot complement the only computed irreducible $SL(4,\C)$-representations are those coming from $\rho_4\circ\iota, \overline{\rho_4}\circ\iota$ and $(\rho_2\otimes\overline{\rho_2})\circ\iota$. (There are however more $PSL(4,\C)$-representations, see \cite[Example 10.2]{gtz}.)
So in this case one should actually have no more than 3 components in the $SL(4,\R)$-character varieties of the figure eight knot complement. For general knots, however, there will be more than 3 components and in \hyperref[compo]{Section \ref*{compo}} we discuss some methods for constructing some of them.

Concerning $SL(n,\C)$-character varieties we use a construction of Ohtsuki-Riley-Sakuma to show the existence of knot complements having arbitrarily many components with vanishing Chern-Simons invariant in their $SL(n,\C)$-character variety.
\begin{thm}\label{arbi} For any natural numbers $N$ and $m$ there exist 2-bridge knots whose $SL(m,\C)$-character variety has more than $N$ connected components and such that the Chern-Simons invariant vanishes on $N$ different components.\end{thm}

The paper is organised as follows. \hyperref[reco]{Section \ref*{reco}} recollects known facts, especially the results from \cite{gtz}. \hyperref[compu]{Section \ref*{compu}} is the heart of the paper, it computes the fundamental class (\hyperref[thm1]{Theorem \ref*{thm1}}) 
for representations of the form $\rho\otimes\overline{\rho}$ with $\rho\colon\pi_1M\to SL(2,\C)$. In \hyperref[compo]{Section \ref*{compo}} we discuss some facts and conjectures about $SL(4,\C)$- and $SL(4,\R)$-character varieties and in particular the proof of \hyperref[arbi]{Theorem \ref*{arbi}}. 

I thank Matthias G\"orner and Sebastian Goette for answering some questions about \cite{pto} and \cite{gz}, respectively, and Neil Hoffman for contributing the proof of \hyperref[2brcs]{Proposition \ref*{2brcs}} on mathoverflow. The computations in the paper have been done with the help of the Sage mathematical software (\cite{sage}).

\section{Recollections}\label{reco}
This section recollects the definition of the fundamental class of a representation and also the definition and properties of 
Cheeger-Chern-Simons invariants. We are going to use the recent approach of Garoufalidis-D.Thurston-Zickert 
which (in the case of flat bundles over triangulated $3$-manifolds) gives a practical 
approach to the computation of a fundamental class of the holonomy representation as an element in the extended Bloch group $\widehat{\mathcal B}(\C)$. (The Cheeger-Chern-Simons invariant of the flat bundle is then obtained by applying the extended Rogers' dilogarithm to that element of $\widehat{\mathcal B}(\C)$.) We describe their approach in \hyperref[ccsgtz]{Section \ref*{ccsgtz}} in some detail because our proof of \hyperref[thm1]{Theorem \ref*{thm1}} depends on this construction. 

In \hyperref[defccs]{Section \ref*{defccs}} we state the general definition of CCS-invariant according to Cheeger and Simons and we explain in \hyperref[invpol]{Lemma \ref*{invpol}} the (apparently well-known) fact that their pull-backs under the geometric representations $\rho_n$ and their tensor products can be computed by using invariant polynomials.

In \hyperref[properties]{Section \ref*{properties}} we discuss some properties of the CCS-invariants (which follow easily from \cite{cs} and \cite{gtz} and are likely to be well-known). 

In \hyperref[symmetries]{Section \ref*{symmetries}} we derive some symmetries in the extended pre-Bloch group (reflecting symmetries of the extended Rogers' dilogarithm) which do not appear in the literature, but which will turn out to be crucial for the proof of \hyperref[thm1]{Theorem \ref*{thm1}}. 

Finally, \hyperref[figureeight]{Section \ref*{figureeight}} discusses the example of the figure eight knot complement and describes in this computationaly easy case the computations that in generality will be performed in the following \hyperref[compu]{Section \ref*{compu}}. 

{\em Conventions}: Throughout the paper $\log(z)$ will mean the branch of the logarithm of $z\in\C$ with imaginary part  $$-\pi< Im(\log(z)) \le\pi.$$ For a manifold $M$ we will always assume to have fixed a basepoint $m_0\in M$ and a lift $\widetilde{m}_0\in\widetilde{M}$, and hence an action of $\pi_1M:=\pi_1(M,m_0)$ on the universal covering space $\widetilde{M}$.

\subsection{Definition of the fundamental class of a representation}\label{ccsgtz}
The group $H_3(SL(2,\C);\Z)$ has an explicit description (by the work of Neumann) as the so-called extended 
Bloch group\footnotemark\footnotetext{Following \cite{dz} and \cite{gz} we define the extended pre-Bloch 
group as what Neumann in \cite{neu} calls the more extended pre-Bloch group.} $\widehat{\mathcal{B}}(\C)$, 
which is a certain subgroup of the extended pre-Bloch group $\PC(\C)$. Together with the Suslin-Sah 
isomorphism $$H_3(SL(n,\C);\Z)\cong H_3(SL(2,\C);\Z)\oplus K_3^M(\C)\ \forall n\ge 3$$
one obtains for $n\ge 3$ a decomposition 
$$H_3(SL(n,\C);\Z)\cong \widehat{\mathcal{B}}(\C)\oplus K_3^M(\C)$$

It is a classical fact (probably first appearing in the work of Dupont-Sah) that the volume of a hyperbolic $3$-manifold can be computed by applying the Bloch-Wigner dilogarithm to a fundamental class in the Bloch group ${\mathcal B}(\C)$. However the Chern-Simons invariant can not be computed from that fundamental class and this is one reason why Neumann in \cite{neu} defined the extended Bloch group $\widehat{\mathcal B}(\C)$ and showed that the volume and Chern-Simons invariant can be computed by applying the extended Rogers' dilogarithm to a fundamental class in $\widehat{\mathcal B}(\C)$. In \cite{gtz}, Garoufalidis, D.Thurston and Zickert extended this approach to $SL(n,\C)$-representations of $3$-manifold groups, i.e., to such a representation they associated a fundamental class in $\widehat{\mathcal B}(\C)$ such that again application of the Rogers' extended dilogarithm computes volume and Chern-Simons invariant of the representation (that is, of the flat bundle whose holonomy is that representation). In this subsection we will review the definitions and results from \cite{gtz}.

\begin{df}\label{extprebloch}The extended pre-Bloch group $\PC(\C)$
is the free abelian group on the set $$\widehat{\C}=\left\{(e,f)\in\C^2\colon exp(e)+exp(f)=1\right\}$$ 
modulo the relations $$
(e_0,f_0)-(e_1,f_1)+(e_2,f_2)-(e_3,f_3)+(e_4,f_4)=0$$
whenever the equations $$
e_2=e_1-e_0,$$
$$e_3=e_1-e_0-f_1+f_0,f_3=f_2-f_1,$$
$$
e_4=f_0-f_1,f_4=f_2-f_1+e_0$$
hold.\end{df}

One should pay attention that $\left[z;2p,2q\right]$ in the notation of \cite{neu}, \cite{dz} and \cite{gz} corresponds to $$(e,f)=(\log(z)+2p\pi i,\log(1-z)-2q\pi i)$$
and hence to $\left[z;2p,-2q\right]$ in the notation of \cite{gtz}. (Here $z\in \C\setminus \left\{0,1\right\}$ and $p,q\in \Z$.)

\begin{df}\label{regu}The extended Rogers' dilogarithm $$R\colon\PC(\C)
\to \C/4\pi^2\Z$$ is defined
on generators of $\PC(\C)$ by $$R((\log(z)+2p\pi i,\log(1-z)+2q\pi i)):= Li_2(z)+\frac{1}{2}(\log(z)+2p\pi i)(\log(1-z)-2q\pi i)-\frac{\pi^2}{6},$$
where $Li_2(z)$ denotes the classical dilogarithm.\end{df}
It is proved in \cite {neu} and \cite{gz} that $R$ is well-defined and a homomorphism. The relation to the Rogers' dilogarithm ${\mathcal R}$ is given by the equality 
$$R((\log(z)+2p\pi i,\log(1-z)+2q\pi i))={\mathcal R}(z)+(p\log(1-z)-q\log(z))\pi i-2pq\pi^2.$$

\begin{df}A closed $3$-cycle is a cell complex $K$ obtained from a finite collection of ordered $3$-simplices by order preserving simplicial gluing maps between pairs of faces. We call $K$ a generalized ideal triangulation of a compact $3$-manifold $M$ if it is homeomorphic to the space $\widehat{M}$ obtained from $M$ by collapsing each boundary component to one point.\end{df}

\begin{df}Let $G$ be a Lie group, $N\subset G$
 a subgroup and $M$ a $3$-manifold. A $(G,N)$-representation is a representation $\pi_1M\to G$ which sends each peripheral subgroup to a conjugate of $N$.\end{df}
 
\begin{df}\label{l}Let $M$ be a compact $3$-manifold (possibly with boundary) and $K$ a generalized ideal triangulation. Then $L$ denotes the corresponding generalized triangulation of the space obtained from the universal covering $\widetilde{M}$
by collapsing each boundary component of $\widetilde{M}$ to a point. \end{df}
We remark that the action of $\pi_1M$ on $\widetilde{M}$ extends to an action of $\pi_1M$ on $L$.

\begin{df}\label{deco}Let $M$ be a compact $3$-manifold (possibly 
with boundary)
and
$K$ a generalized ideal triangulation.
A decoration of a $(G,N)$-representation $\rho\colon
\pi_1M\to G$ is a $\rho$-equivariant assignment 
$$L_0\to G/N,$$
i.e.,
associating an $N$-coset to each vertex of $L$ such that if $\alpha\in \pi_1M$ and the coset $g_vN$
is associated to $v$, then the coset 
$\rho(\alpha)g_vN$
is associated to $\alpha v$.

We will say that a simplex $\sigma=(v_0,v_1,v_2,v_3)$
from $K$ is decorated by the tuple $(g_0N,g_1N,g_2N,g_3N)$.
\end{df}

In what follows $G$ will be a subgroup of $GL(n,\C)$ and $N\subset G$ will be the subgroup of upper triangular matrices with all diagonal entries equal to $1$.

\begin{df}\label{ptolemy1} Let $K$ be a generalized ideal triangulation of a compact $3$-manifold $M$ and let $\rho\colon \pi_1M\to G$ be a decorated $(G,N)$-representation.

The Ptolemy coordinates $\left\{c_t\right\}_t$ of a decorated $3$-simplex $$(g_0N,g_1N,g_2N,g_3N)$$ 
are the assignment
$$t=(t_0,t_1,t_2,t_3)
\to c_t:=det (\bigcup_{i=0}^3\left\{
g_i\right\}_{t_i})$$
for each $4$-tuple $t=(t_0,t_1,t_2,t_3)$ of nonnegative integers with $$t_0+t_1+t_2+t_3=n.$$

Here $\left\{g_i\right\}_{t_i}$ means the (ordered) set of the first $t_i$ column vectors of $g_i\in GL(n,\C)$ and $\bigcup_{i=0}^3\left\{
g_i\right\}_{t_i}$
means the matrix whose (ordered) column set is composed by the 
$\left\{g_i\right\}_{t_i}$.\end{df}

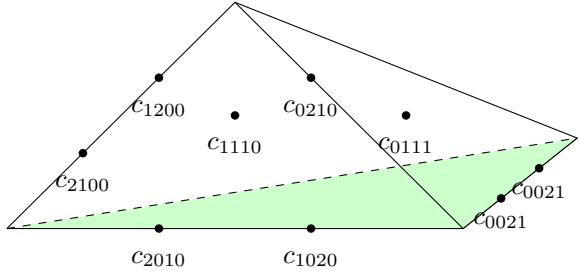
\begin{tikzpicture}
\path[draw](3,3)--(7.5,1.2)--(6,0);
\path[draw, fill=green!20](0,0)--(6,0)--(7.5,1.2);
\path[draw](0,0)--(3,3)--(6,0);
\draw[dashed](0,0)--(7.5,1.2);
\fill (1,1) circle (1.6pt) node [below=5pt] {$c_{2100}$};
\fill (2,2) circle (1.6pt) node [below=5pt] {$c_{1200}$};
\fill (3,1.5) circle (1.6pt) node [below=5pt] {$c_{1110}$};
\fill (2,0) circle (1.6pt) node [below=5pt] {$c_{2010}$};\fill (4,0) circle (1.6pt) node [below=5pt] {$c_{1020}$};
\fill (5.25,1.5) circle (1.6pt) node [below=5pt] {$c_{0111}$};
\fill (6.5,0.4) circle (1.6pt) node [below=2pt] {$c_{0021}$};\fill (7,0.8) circle (1.6pt) node [below=2pt] {$c_{0021}$};
\fill (4,2) circle (1.6pt) node [below=5pt] {$c_{0210}$};
\end{tikzpicture}

One visualizes the Ptolemy coordinates (of a simplex)
by fixing some identification of the simplex with 
$$
\Delta^3_n:=
\left\{
(x_0,x_1,x_2,x_3)\in \R^4
\colon x_i\ge 0, x_0+x_1+x_2+x_3=n
\right\}$$
and by
attaching $c_t$ to the point $(t_0,t_1,t_2,t_3)$. The picture above shows $\Delta^3_3$ with some of the $c_t$ attached.

\begin{df}\label{generi}
A decoration is called generic if all Ptolemy coordinates are nonzero:
$$c_t\not=0\ \forall t\in\Z_{\ge0},t_0+t_1+t_2+t_3=n.$$
\end{df}
One can always obtain generic decorations by performing a barycentric subdivision on simplices with nongeneric decorations.

Next we describe some $3$-simplices embedded in $\Delta^3_n$.
For each 
$$
\alpha\in\Delta^3_{n-2}\cap \Z^{4}$$
one can consider the $3$-simplex with vertices corresponding to
$$
\alpha+e_i, i=0,1,2,3.$$
Each of its edges
$$
(\alpha+e_i,\alpha+e_j)$$
has one Ptolemy coordinate attached to it, namely $$
c_{\alpha_{ij}}:=
c_{\alpha+e_i+e_j}.$$

\begin{df}\label{invariant}Let $K
=\cup_{k=1}^r T_k$ be
a generalized ideal triangulation of a compact, orientable
$3$-manifold $M$ and let $\rho\colon \pi_1M\to G$ be a generic 
decorated $(G,N)$-representation, with Ptolemy coordinates $c_t^k
$ for each simplex $T_k$.

For each simplex $T$ and each $\alpha\in
\Delta^3_{n-2}\cap \Z^4$ define 
$$
\setlength{\multlinegap}{0pt}
\begin{array}{c}
\hspace{-1.5in}\lambda(c_\alpha)
=
(\log(c_{\alpha_{03}})+\log(c_{\alpha_{12}})-\log(c_{\alpha_{02}})-\log(c_{\alpha_{13}}),\\
\hspace{1.5in}\log(c_{\alpha_{01}})+\log(c_{\alpha_{23}})-\log(c_{\alpha_{02}})-\log(c_{\alpha_{13}}))\in\PC
(\C).
\end{array}
$$

Then define
$$\lambda(K,\rho)=\sum_{k=1}^r \epsilon_k \sum_{\alpha\in \Delta^3_{n-2}\cap \Z^4}\lambda(c^k_\alpha)\in\PC(\C),$$
where $\epsilon_k$ is $\pm 1$ according to whether the orientation of $T_k$ agrees with the orientation of $M$ or not.\end{df}

It is known (but it will actually play no role for our calculations) that $$\lambda(K,\rho)\in\widehat{\mathcal{B}}(\C):=ker(\hat{\nu})$$ for the extended Dehn invariant $\hat{\nu}(e,f):=e\wedge f\in\C\wedge_{\Z}\C$. The element $\lambda(K,\rho)$ does not depend on the triangulation and we will denote it by
$$\rho_*\left[M,\partial M\right]\in \widehat{\mathcal{B}}(\C).$$
In fact for $G=SL(2,\C)$ it corresponds to the image of the fundamental class under the epimorphism $H_3(SL(2,\C)^\delta,N^\delta;\Z)\to H_3(SL(2,\C)^\delta;\Z)\cong \widehat{\mathcal{B}}(\C)$.
 
The following result from \cite{gtz} will be useful to avoid too many case distinctions in our arguments.
\begin{pro}\label{lift}(\cite[Proposition 7.7]{gtz}) Under the assumptions of 
\hyperref[invariant]{Definition \ref*{invariant}} let $c$ be the Ptolemy coordinates of the generic decorated representation on $K$. For any lift $\tilde{c}$ of $c$ we have $$\lambda(K,\rho)=\sum_{k=1}^r \epsilon_k \sum_{\alpha\in \Delta^3_{n-2}\cap \Z^4}\tilde{\lambda}(c^k_\alpha)\in\PC(\C).$$\end{pro}
Here, a lift $\tilde{c}$ of $c$ means a choice of logarithm for each $c_\alpha^k$ (i.e., of a complex number whose difference with $\log(c_\alpha^k)$ is an integer multiple of $2\pi i$) such that the choices agree whenever coordinates correspond to glued faces in $K$, and $\tilde{\lambda}(c_\alpha)$ is then defined as $$\tilde{\lambda}(c_\alpha)=(\tilde{c}_{\alpha_{03}}+\tilde{c}_{\alpha_{12}}-\tilde{c}_{\alpha_{02}}-\tilde{c}_{\alpha_{13}},\tilde{c}_{\alpha_{01}}+\tilde{c}_{\alpha_{23}}-\tilde{c}_{\alpha_{02}}-\tilde{c}_{\alpha_{13}})\in\PC(\C).$$
So, in the formula of \hyperref[invariant]{Definition \ref*{invariant}} one can replace $\log$ by any choice of logarithm as long as we make the same choice on common faces or edges of different simplices.

{\bf $pSL(n,\C)$-representations (\cite[Section 6.3]{gtz}).} For boundary-unipotent representations to 
$$pSL(n,\C)=SL(n,\C)/\left\{\pm1\right\}$$ 
the ptolemy coordinates are only well-defined as elements of 
$\C^*/\left\{\pm1\right\}$ and thus $\lambda(c_\alpha)$ takes value in $\PC(\C)_{PSL}$, the free abelian group over 
$$\left\{(e,f)\in\C^2\colon \pm exp(e)\pm exp(f)=1\right\}$$ 
modulo the 5-term relation. 

\subsection{Definition of CCS-invariants}\label{defccs}
For a flat complex vector bundle $\mathcal{V}:E\to X$
Cheeger-Simons (\cite[Section 4]{cs}) define\footnotemark\footnotetext{In their normalization the Chern character is an element of $H^{2k-1}(X;\C/\Z)$.} Chern characters $\hat{c}_k({\mathcal{V}})\in H^{2k-1}(X;\C/4\pi^2\Z)$. In this paper we will be interested in the character $\hat{c}_2(\mathcal{V})$ for flat $SL(n,\C)$-bundles over $3$-manifolds. For a closed, orientable 3-manifold $M$ and a representation $\rho:\pi_1M\rightarrow SL(n,C)$ we consider its Cheeger-Chern-Simons invariant $$CCS(M,\rho)=\int_M \hat{c}_2(\mathcal{V}_\rho),$$
where ${\mathcal{V}}_\rho$ means the flat $n$-dimensional complex vector bundle over $M$ with holonomy $\rho$. An explicit formula is 
$$CCS(M,\rho)=\frac{1}{2}\int_Ms^*(Tr(\theta\wedge d\theta+\frac{2}{3}\theta\wedge\theta\wedge\theta))\ mod\ 4\pi^2\Z,$$
where $\theta$ is a flat connection and $s$ a section of $\mathcal{V}_\rho$, which exists because $SL(n,\C)$ is 2-connected.

The universal Cheeger-Chern-Simons class $\hat{c}_2$ of flat $SL(n,\C)$-bundles 
is defined in \cite{cs} as an element in $H^3(SL(n,\C);\C/4\pi^2\Z)$. 
To get a Cheeger-Chern-Simons invariant also for cusped manifolds one lets $N\subset SL(n,\C)$ be the subgroup of upper triangular matrices with 1's on the diagonals and 
uses the map $ H^3(SL(n,\C),N;\C/4\pi^2\Z)\to H^3(SL(n,\C);\C/4\pi^2\Z)$ to consider $\hat{c}_2$ as a relative class $$\hat{c}_2\in H^3(BSL(n,\C)^\delta, BN^\delta;\C/4\pi^2\Z),$$
see \cite[Section 6.1]{gtz}. Then for a flat bundle with boundary-unipotent holonomy (i.e.\ the restriction of the holonomy to $\partial M$ having image in conjugates of $N$) one defines the Cheeger-Chern-Simons invariant via the pullback\footnotemark\footnotetext{Although the induced homomorphism $\pi_1M\to SL(n,\C)$ sends the fundamental groups of different boundary components to possibly different conjugates of $N$, nonetheles one has a well-defined homomorphism $H^3(BSL(n,\C)^\delta,BN^\delta;\C/4\pi^2\Z)\to H^3(M,\partial M;\C/4\pi^2\Z)$.} of $\hat{c}_2$ under the classifying map $M\to BSL(n,\C)^\delta$.

As $\C/4\pi^2Z$ is divisible, one may consider $\hat{c}_2$ as a homomorphism
$$\hat{c}_2\colon H_3(SL(n,\C);\Z)\to \C/4\pi^2Z.$$
The group $H_3(SL(2,\C);\Z)$ has an explicit description (by the work of Neumann) as the so-called extended Bloch 
group\footnotemark\footnotetext{Following \cite{dz} and \cite{gz} we define the extended pre-Bloch group as what 
Neumann in \cite{neu} calls the more extended pre-Bloch group.} $\widehat{\mathcal{B}}(\C)$, which is a certain 
subgroup of the extended pre-Bloch group $\PC(\C)$ described in \hyperref[ccsgtz]{Section \ref*{ccsgtz}}. 
Together with the Suslin-Sah isomorphism $H_3(SL(n,\C);\Z)\cong H_3(SL(2,\C);\Z)\oplus K_3^M(\C)$ (for $n\ge 3$) one obtains a decomposition 
$$H_3(SL(n,\C);\Z)\cong \widehat{\mathcal{B}}(\C)\oplus K_3^M(\C)$$
and it turns out that $\hat{c}_2$ vanishes on the Milnor K-theory $K_3^M(\C)$ (see \cite[Section 8]{gtz}), thus $\hat{c}_2$ depends only on its values on $\widehat{\mathcal{B}}(\C)$. As explained in \hyperref[ccsgtz]{Section \ref*{ccsgtz}}, Garoufalidis, D.\ Thurston and Zickert associate to each flat $SL(n,\C)$-bundle over a $3$-manifold (with unipotent holonomy at the boundary) a "fundamental class" 
$$(B\rho)_*\left[M,\partial M\right]\in\widehat{\mathcal{B}}(\C)$$ 
and they exhibit an explicit method for computing $\hat{c}_2$ on this element, 
see \hyperref[gtz]{Proposition \ref*{gtz}} below. This yields a computable formula for the Cheeger-Chern-Simons invariant of the associated flat bundle. 

Using the Cheeger-Chern-Simons invariant, the volume and Chern-Simons invariant of a representation are defined as follows.
\begin{df}\label{ccsdef} (\cite[Definition 2.11]{gtz}): For a compact, orientable, aspherical $3$-manifold $M$ and a boundary-unipotent representation $\rho\colon\pi_1M\to SL(n,\C)$ one defines the volume\footnotemark\footnotetext{This is not the same as the volume of representations defined via pulling back the volume form of the symmetric spaces.} and Chern-Simons invariant of $\rho$ by 
$$-CS(\rho)+i Vol(\rho)=\langle \hat{c}_2,(B\rho)_*\left[M,\partial M\right]\rangle,$$
where $$(B\rho)_*\colon H_3(M,\partial M;\Z)\to H_3(BSL(n,\C)^\delta,BN^\delta;\Z)$$ is induced by the classifying map $B\rho\colon M\to BSL(n,\C)^\delta$ of $\rho$ (i.e., of the associated flat bundle).\end{df}
The motivation for this naming is Yoshida's theorem (see \cite[Theorem 2.8]{gtz}) which implies that for a closed hyperbolic $3$-manifold and a lift $\iota\colon\pi_1M\to SL(2,\C)$ of its geometric representation, $Vol(\iota)$ (as defined above) is the hyperbolic volume $Vol(M)$ and $CS(\iota)$ is the Chern-Simons invariant $CS(M)$ of the Levi-Civita connection for the hyperbolic metric. (The analogous result for cusped hyperbolic $3$-manifolds is true modulo $\pi^2$ and is proved in \cite[Corollary 14.6]{neu}.)

The following lemma seems well-known, although we were not able to locate a reference in the literature. For the proof we will need the universal approach 
to CCS-classes from \cite{cs}, which is as follows. For a Lie 
group $G$, let $I^*(G)$ be its invariant polynomials over $\C$, $BG$ its classifying 
space, $CW_k\colon I^k(G)\to H^{2k}(BG;\C)$ the Chern-Weil isomorphism 
and $$K^{2k}(G;\Z)=\left\{(P,u)\in I^k(G)\times H^{2k}(BG;\Z)\colon CW_k(P)=u_{\C}\right\}.$$ 
Then, for each principal 
$G$-bundle with flat connection $\theta$ and base $B$, there is 
a natural map $$S_{(P,u)}(\theta)\colon K^{2k}(G;\Z)\to H^{2k-1}(B;\C/\Z).$$ 
The Cheeger-Chern-Simons class $\hat{c}_2$ of the flat bundle equals 
$4\pi^2S_{P,u}(\theta)$ with $P\in I^2(SL(n,\C))$ the invariant polynomial 
$P(A)=-\frac{1}{4\pi^2}Tr(A^2)$ and $u\in H^4(BSL(n,\C);\Z)$ the universal second Chern class.
\begin{lem}\label{invpol}
Let $M$ be a compact, orientable $3$-manifold and $\rho\colon\pi_1M\to SL(2,\C)$ any representation. For $n\ge 2$ let $\rho_n\colon SL(2,\C)\to SL(n,\C)$ be the irreducible representation corresponding to the unique $n$-dimensional $\C$-linear of the Lie algebra $\mathfrak{sl}(2,\C)$. Then \\
a) $$CCS(M,\rho_n\circ\rho)=
\left(\begin{array}{c}n+1\\
3\end{array}\right) CCS(M,\rho),$$
b) $$CCS(M,(\rho_n\otimes\overline{\rho_m})\circ\rho)=mCCS(M,\rho)+n\overline{CCS(M,\rho)}.$$
\end{lem}
\begin{pf}
By the previous paragraph $CCS$ is defined via the invariant polynomial 
$P(A)=-tr(A^2)\in I^2(SL(n,\C))$.
Naturality of $S_{(P,u)}(\theta)$ means $$S_{(\rho_n^*P,\rho_n^*u)}(\theta)=S_{(P,u)}(\pi_{n}\theta)$$
(see \cite[Theorem 2.2]{cs}), 
when $$\pi_n=D_e\rho_n\colon \mathfrak{sl}(2,\C)\to\mathfrak{sl}(n,\C)$$ is the induced homomorphism of Lie algebras, and it implies that $CCS(M,\rho_n\circ\rho)$ can be computed by integrating the invariant polynomial $\rho_n^*P$ instead of $P$.

In view of naturality and of $\rho_n^*P(A)=P(\pi_n(A))$ we then just 
have to show that $$P(\pi_n(A))=\left(\begin{array}{c}n+1\\
3\end{array}\right)P(A)$$ and 
$$P(\pi_n(A)\otimes 1_m+1_n\otimes \overline{\pi_m(A)})=mP(A)+n\overline{P(A)},$$ 
where we have used that the Lie algebra homomorphism induced from the Lie group homomorphism $B\to B\otimes B$ is $A\to A\otimes 1+1\otimes A$.

In the first equality we have a $\C$-linear map, so it suffices to check the equality on the $\C$-basis of $\mathfrak{sl}(2,\C)$ given by $X=\left(\begin{array}{cc}0&1\\
0&0\end{array}\right), Y=\left(\begin{array}{cc}0&0\\
1&0\end{array}\right), H=\left(\begin{array}{cc}1&0\\
0&-1\end{array}\right)$. One easily checks that $\pi_n(X)$ 
is an upper triangular matrix and $\pi_n(Y)$ is a lower triangular 
matrix (both with $0$'s on the diagonal), so we have $P(\pi_n(X))=0$ and $P(\pi_n(Y))=0$ as desired. 
Moreover, $\pi_n(H)$ is the diagonal matrix $diag(n-1,n-3,\ldots,-n+1)$, from which one obtains $Tr(\pi_n(H))^2=2\left(\begin{array}{c}n+1\\
3\end{array}\right).$

For the second claim it suffices by symmetry, additivity and compatibility with complex conjugation to show that $P(\pi_n(A)\otimes 1_m)=mP(A)$ and again to consider the $\C$-basis. For $A=X$ and $A=Y$ one again obtains triangular matrices with $0$'s on the diagonal, so $P(\pi_n(X)\otimes 1_m)=0$ and $P(\pi_n(Y)\otimes 1_m)=0$. For $A=H$ one observes that $\pi_n(H)\otimes 1_m$ has the same eigenvalues as $\pi_n(H)$ but each occuring with multiplicity $m$, so $$Tr(\pi_n(H)^2\otimes 1_m)=m Tr(\pi_n(H)^2)$$
which yields the claim.

\end{pf}

The previous lemma suggests the following conjecture.
\begin{conj}\label{tensconj}Under the assumptions of \hyperref[invpol]{Lemma \ref*{invpol}} we have
$$(\rho_n\otimes\overline{\rho_m})_*\left[M,\partial M\right]=m\left[M,\partial M\right]+n\overline{\left[M,\partial M\right]}.$$\end{conj}
In this paper we will handle the case $n=m=2$.

\subsection{Properties of CCS-invariants}\label{properties}

In this subsection we recollect some properties of the CCS-invariant.
The following is the main result of Garoufalidis-D.Thurston-Zickert in \cite{gtz}. It shows that the Cheeger-Chern-Simons invariant (see \hyperref[ccsdef]{Definition \ref*{ccsdef}}) can be computed from $\rho_*\left[M,\partial M\right]$. 
\begin{pro}\label{gtz}(\cite[Theorem 1.3]{gtz}) Let $K=\cup_{k=1}^r T_k$ be
a generalized ideal triangulation of a compact, orientable
$3$-manifold $M$ (possibly with boundary) and let $\rho\colon \pi_1M\to G$ be a generic 
decorated $(SL(n,\C),N)$-representation, with ptolemy coordinates $c_t^i
$ for each simplex $T_i$. Then
$$R(\lambda(K,\rho))=-CS(\rho)+i Vol(\rho)\in\C/4\pi^2\Z,$$
where $R$ denotes the extended Rogers' dilogarithm from \hyperref[regu]{Definition \ref*{regu}} and $\lambda(K,\rho)\in \PC(\C)$ was defined in \hyperref[invariant]{Definition \ref*{invariant}}.\end{pro}
For boundary-unipotent representations to $pSL(n,\C)$ we have the fundamental class in $\PC(\C)_{PSL}$ (see the last paragraph of \hyperref[ccsgtz]{Section \ref*{ccsgtz}}) and the extended Rogers' dilogarithm is then well-defined modulo $\pi^2\Z$. The ptolemy coordinates in \hyperref[ptolemy]{Definition \ref*{ptolemy}} are defined as elements of $\C^*/\left\{\pm 1\right\}$, so $\lambda(K,\rho)$ in \hyperref[invariant]{Definition \ref*{invariant}} can (by taking imaginary parts of logarithms between $0$ and $\pi$) be defined as an element in $\PC(\C)_{PSL}$ and \hyperref[gtz]{Proposition \ref*{gtz}} holds as an equality in $\C/\pi^2\Z$, see \cite[Section 9.2]{gtz}.

While the Cheeger-Chern-Simons invariant is in general not additive for direct sums of $GL(n,\C)$-bundles, additivity holds for $SL(n,\C)$-bundles. 
\begin{lem}\label{additivity}Let $\rho_1,\rho_2$ be boundary-unipotent representations from $\pi_1M$ to $SL(n,\C)$, for a compact manifold $M$. Then$$CCS(M,\rho_1\oplus\rho_2)=CCS(M,\rho_1)+CCS(M,\rho_2)$$
\end{lem}
\begin{pf}
From Cheeger-Simons (\cite[Theorem 4.6]{cs}) follows
$$\hat{c_2}({\mathcal{V}}\oplus{\mathcal{W}})=\hat{c_2}({\mathcal{V}})+\hat{c_2}({\mathcal{W}})+\hat{c_1}({\mathcal{V}})*\hat{c_1}({\mathcal{W}})$$
for a certain multiplication $*$ defined in \cite[page 56]{cs}.

It is well-known that $\hat{c}_1$ vanishes for all flat $SL(n,\C)$-bundles. Indeed, $$\hat{c}_1\in H^1(BGL(n,\C);\C/2\pi i\Z)$$ 
is represented by the cocycle $g\to \log(det(g))$. The claim follows.\end{pf}\\

A direct consequence of the explicit formula in \cite{gtz} is the compatibility of $CCS$ with complex conjugation.
\begin{lem}\label{conj}For any boundary-unipotent representation $\rho\colon\pi_1M\to SL(n,\C)$ we have $$CCS(M,\overline{\rho})=\overline{CCS(M,\rho)}.$$\end{lem}
\begin{pf}Going through the formulas in \cite{gtz} one sees that a decoration for $\overline{\rho}$ can be obtained by 
applying complex conjugation to a decoration for $\rho$, and that the Ptolemy coordinates of these decorations are 
related by complex conjugation. According to \hyperref[lift]{Proposition \ref*{lift}} the value of $CCS(M,\rho)$ does 
not depend on the lifts of the Ptolemy coordinates, so one may choose the lifts of the Ptolemy coordinates 
for $\overline{\rho}$ to be exactly the complex conjugates of the lifts of the Ptolemy coordinates for $\rho$. (For example one may choose $\tilde{c}=\log(c)$ whenever $c\not\in \R_{<0}$, and for all $c\in\R_{<0}$ one may choose $\tilde{c}=\log(c)$ for $\rho$, but $\tilde{c}=\log(c)-2\pi i=\overline{\log(c)}$  for $\overline{\rho}$.) 

The formula in \hyperref[invariant]{Definition \ref*{invariant}} then implies $\lambda(K,\overline{\rho})=\overline{\lambda(K,\rho)}$ and now the claim follows from 
\hyperref[gtz]{Proposition \ref*{gtz}} and the equality $R(\overline{e},\overline{f})=\overline{R(e,f)}$ which is immediate from \hyperref[regu]{Definition \ref*{regu}}.\end{pf}\\

An immediate consequence is that a boundary-unipotent representation which can be conjugated to a representation in $SL(n,\R)$ must have vanishing volume. It is perhaps worth-mentioning that the ptolemy coordinates of such a representation are not necessarily real, basically because the peripheral subgroups are conjugate to $N\cap SL(n,\R)$ inside $SL(n,\C)$ but not necessarily inside $SL(n,\R)$. The computations in \hyperref[compu]{Section \ref*{compu}} actually provide an example of this phenomenon.

For the proof of \hyperref[2br]{Proposition \ref*{2br}} and hence \hyperref[arbi]{Theorem \ref*{arbi}} we will also need the following straightforward equality.

\begin{lem}\label{degree}If $f\colon(M_1,\partial M_1)\to (M_2,\partial M_2)$ has mapping degree $deg(f)$, then $$CCS(M_2,\rho)=\frac{1}{deg(f)}CCS(M_1,\rho\circ f_*)$$
for any boundary-unipotent representation $\rho\colon\pi_1M_2\to SL(n,\C)$ and the induced homomorphism $f_*\colon\pi_1M_1\to\pi_1M_2$.\end{lem}
\begin{pf} This is immediate from $$CCS(M_2,\rho)=\langle\hat{c}_2,(B\rho)_*\left[M_2,\partial M_2\right]\rangle$$ and $$\left[M_2,\partial M_2\right]=deg(f)\left[M_1,\partial M_1\right].$$\end{pf}

\subsection{Symmetries of the extended Rogers' dilogarithm}\label{symmetries}
The aim of this subsection is to prove some relations in the extended pre-Bloch group $\PC(\C)$, which are more complicated than the well-known relations in the pre-Bloch group ${\mathcal P}(\C)$. These relations will play a central role in the proof of \hyperref[thm1]{Theorem \ref*{thm1}} in \hyperref[compu]{Section \ref*{compu}}.

\begin{df}For $e\in\C$ define $\chi(e)\in\PC(\C)$ by $$\chi(e)=(e,f+2\pi i)-(e,f)$$
where $f\in\C$ is some complex number satisfying $exp(e)+exp(f)=1$.\end{df}
In the following lemma, $R\colon\PC(\C)\to\C/4\pi^2\Z$ denotes the extended Rogers' dilogarithm from \hyperref[regu]{Definition \ref*{regu}}.
\begin{lem}\label{hom}
i) $\chi$ is a homomorphism $ \C/4\pi i\Z\to \PC(\C)$
with respect to the additive structures on $\C$ and $\PC(\C)$.

ii) $R(\chi(e))=-\pi i e\ mod\ 4\pi^2$ for all $e\in \C$. 

iii) $R$ is injective on the image of $\chi$, and $im(\chi)=ker(\PC(\C)\to\mathcal{P}(\C))$.

iv) For $(e,f)\in\PC(\C)$ and $p,q\in\Z$ we have 
$$(e+2p\pi i,f+2q\pi i)-(e,f)=\chi(2pq\pi i+qe-pf).$$


\end{lem}
\begin{pf}
i)-iii) are (in a slightly different language) proved in \cite[Theorem 3.12]{gz}. (There is a different sign in ii) because of the notational difference explained in the remark before \hyperref[regu]{Definition \ref*{regu}}.) Equation iv) follows from i)-iii), for general fields $F$ it is also proved as a consequence of the 5-term relation in \cite[Lemma 3.16]{zic}. 
\end{pf}\\

The following relations will be crucial for the computations in \hyperref[compu]{Section \ref*{compu}}, in particular they will be used in the proof of \hyperref[cancel]{Corollary \ref*{cancel}} and thus \hyperref[thm1]{Theorem \ref*{thm1}}.

\begin{lem}\label{gz} The following relations hold whenever $Im(z)>0$ and $p,q\in\Z$

i)
$(\log(\frac{1}{z})-2p\pi i,\log(1-\frac{1}{z})+2(q-p)\pi i)$
$$=-(\log(z)+2p\pi i,\log(1-z)+2q\pi i)+\chi(-\frac{1}{2}\log(z)+(2p^2+p)\pi i)$$


ii) 
$(\log(1-z)+2q\pi i,\log(z)+2p\pi i)$
$$=-(\log(z)+2p\pi i,\log(1-z)+2q\pi i)+\chi(-\frac{\pi i}{6})$$

iii) 
$(\log(\frac{1}{1-z})-2q\pi i,\log(\frac{-z}{1-z})+
2(p-q)\pi i)$
$$=(\log(z)+2p\pi i,\log(1-z)+2q\pi i)+
\chi(\frac{1}{2}\log(1-z)+(2q^2-q+\frac{1}{6})\pi i)$$

iv) 
$(\log(1-\frac{1}{z})+
2(q-p)\pi i,\log(\frac{1}{z})-2p\pi i)$
$$=(\log(z)+2p\pi i,\log(1-z)+2q\pi i)+\chi(\frac{1}{2}\log(z)+(2p^2-p-\frac{1}{6})\pi i)$$

v) $(\log(-\frac{\overline{z}}{1-\overline{z}})+2(p-q)\pi i,\log(\frac{1}{1-\overline{z}})-2q\pi i)$
$$=
-(\log(\oz)+2p\pi i,\log(1-\oz)+2q\pi i)+\chi(\frac{1}{2}\log(1-\oz)+(2q^2-q-\frac{1}{3})\pi i$$
\end{lem}
\begin{pf}

i) From \hyperref[regu]{Definition \ref*{regu}} we get $$R(\log(z),\log(1-z))+R(\log(\frac{1}{z}),\log(1-\frac{1}{z}))=$$
$$Li_2(z)+\frac{1}{2}\log(z)\log(1-z)-\frac{\pi^2}{6}+Li_2(\frac{1}{z})+\frac{1}{2}\log(\frac{1}{z})\log(1-\frac{1}{z})-\frac{\pi^2}{6}\ \mbox{mod}\ 4\pi^2.$$

From \cite[Section 2]{zag} we have $$Li_2(z)+Li_2(\frac{1}{z})-\frac{\pi^2}{3}=-\frac{\pi^2}{2}-\frac{1}{2}\log^2(-z).$$
Moreover we have $\log(\frac{1}{z})=-\log(z), \log(1-\frac{1}{z})=\log(z-1)-\log(z)$ because of $\mid arg(z-1)-arg(z)\mid<\pi$, and for $Im(z)>0$ we have $\log(1-z)=\log(z-1)-\pi i$ and $\log(-z)=\log(z)-\pi i$. Thus
the above expression simplifies to
$$-\frac{1}{2}(\log(z)-\pi i)^2-\frac{\pi^2}{2}+\frac{1}{2}\log(z)\log(1-z)+\frac{1}{2}(-\log(z))(\log(1-z)-\log(z)+\pi i)\ \mbox{mod}\ 4\pi^2$$
$$=\frac{1}{2}\pi i\log(z)=R(\chi(-\frac{1}{2}\log(z)
))\ \mbox{mod}\ 4\pi^2,$$
thus $$(\log(z),\log(1-z))+(\log(\frac{1}{z}),\log(1-\frac{1}{z}))=\chi(-\frac{1}{2}\log(z))
.$$

Then
we apply \hyperref[hom]{Lemma \ref*{hom}}, iv) to get
$$(\log(z)+2p\pi i,\log(1-z)+2q\pi i)+(\log(\frac{1}{z})-2p\pi i,\log(1-\frac{1}{z})+(2q-2p)\pi i)=$$
$$\setlength{\multlinegap}{0pt}
\begin{array}{c}
\shoveleft{
(\log(z),\log(1-z))+\chi(2pq\pi i+q\log(z)-p\log(1-z))+(\log(\frac{1}{z}),\log(1-\frac{1}{z}))}\\
\hspace{1.5in}+\chi(2p(q-p)\pi i+(q-p)\log(\frac{1}{z})+p\log(1-\frac{1}{z})\end{array}$$

$$=\chi(-\frac{1}{2}\log(z))+\chi(2p^2\pi i+(q-(q-p))\log(z)-p\log(1-z)+p\log(1-\frac{1}{z}))$$

We are assuming $Im(z)>0$, which implies 
$\log(z-1)-\log(1-z)=\pi i$. Moreover $z$ and $z-1$ have positive imaginary parts, which implies that their arguments differ by less than $\pi$, so $\log(\frac{z-1}{z})=\log(z-1)-\log(z)$. So the above sum simplifies to
$$=\chi(-\frac{1}{2}\log(z)+(2p^2+p)\pi i)$$

ii) We have from \cite[Section 2]{zag} that
$$Li_2(z)+Li_2(1-z)=\frac{\pi^2}{6}-\log(z)\log(1-z),$$
hence
$$R(\log(z),\log(1-z))+R(\log(1-z),\log(z))=\frac{\pi^2}{6}-\log(z)\log(1-z)+2\frac{1}{2}\log(z)\log(1-z)-2\frac{\pi^2}{6}=-\frac{\pi^2}{6},$$
so $$(\log(z),\log(1-z))+(\log(1-z),\log(z))=\chi(-\frac{\pi i}{6}).$$

Then we apply \hyperref[hom]{Lemma \ref*{hom}}, iv) to get
$$(\log(z)+2p\pi i,\log(1-z)+2q\pi i)+(\log(1-z)+2q\pi i,\log(z)+2p\pi i)=$$
$$(\log(z),\log(1-z))+\chi(2pq\pi i+q\log(z)-p\log(1-z))+(\log(1-z),\log(z))+\chi(2pq\pi i
+p\log(1-z)-q\log(z))$$
$$=(\log(z),\log(1-z))+(\log(1-z),\log(z))=\chi(-\frac{\pi i}{6}).$$
because of $\chi(4pq\pi i)=0$. (For equation ii) one
actually does not
need the asumption $Im(z)>0$.)

iii) Because of
$Im(\frac{1}{1-z})>0$ one can apply i) to $\frac{1}{1-z}$ and get
$$(\log(\frac{1}{1-z})-2q\pi i,\log(\frac{-z}{1-z})+
2(p-q)\pi i)$$
$$=-(\log(1-z)+2q\pi i,\log(z)+2p\pi i)+\chi(-\frac{1}{2}\log(\frac{1}{1-z})+(2q^2+q)\pi i)$$
$$=(\log(z)+2p\pi i,\log(1-z)+2q\pi i)+\chi(\frac{1}{2}\log(1-z)+(2q^2-q)\pi i+\frac{\pi i}{6}).
$$

iv) Application of ii) to $1-\frac{1}{z}$ yields
$$(\log(1-\frac{1}{z})+2(q-p)
\pi i,\log(\frac{1}{z})-2p
\pi i)=-(\log(\frac{1}{z})-2p
\pi i,\log(1-\frac{1}{z})+2(q-p)
\pi i)+\chi(-\frac{\pi i}{6})$$
$$=(\log(z)+2p\pi i,\log(1-z)+2q\pi i)+\chi(\frac{1}{2}\log(z)+(2p^2-p)\pi i-\frac{\pi i}{6}).$$

v) Because of $Im(1-\overline{z})>0$ one can apply iv) to $1-\overline{z}$ and use $1-\frac{1}{1-\oz}=-\frac{\oz}{1-\oz}$ to get
$$(\log(-\frac{\overline{z}}{1-\overline{z}})+2(p-q)\pi i,\log(\frac{1}{1-\overline{z}})-2q\pi i)$$
$$=(\log(1-\oz)+2q\pi i,\log(\oz)+2p\pi i)+\chi(\frac{1}{2}\log(1-\oz)+(2q^2-q-\frac{1}{6})\pi i)$$
which by ii) equals to
$$-(\log(\oz)+2p\pi i,\log(1-\oz)+2q\pi i)+\chi(\frac{1}{2}\log(1-\oz)+(2q^2-q-\frac{1}{3})\pi i).$$

\end{pf}

Remark: Other relations have been proved in \cite[Proposition 13.1]{neu} and \cite[Proposition 5.1]{gz}, but they appear not to be correct.
\begin{cor}\label{coro} For $Im(z)<0$ and $p,q\in\C$ we have:\\
i)
$(\log(\frac{1}{z})-2p\pi i,\log(1-\frac{1}{z})+2(q-p)\pi i)$
$$=-(\log(z)+2p\pi i,\log(1-z)+2q\pi i)+\chi(\frac{1}{2}\log(z)-(2p^2+p)\pi i)$$
ii) $(\log(\frac{1}{1-z})-2q\pi i,\log(\frac{-z}{1-z})+
2(p-q)\pi i)$
$$=(\log(z)+2p\pi i,\log(1-z)+2q\pi i)+
\chi(-\frac{1}{2}\log(1-z)-(2q^2-q-\frac{1}{6})\pi i)$$
iii) $(\log(1-\frac{1}{z})+
2(q-p)\pi i,\log(\frac{1}{z})-2p\pi i)$
$$=(\log(z)+2p\pi i,\log(1-z)+2q\pi i)+\chi(\frac{1}{2}\log(z)-(2p^2+p+\frac{1}{6})\pi i)$$
\end{cor}
\begin{pf}We can apply \hyperref[gz]{Lemma \ref*{gz}} to $\frac{1}{z}$ or $1-z$, respectively.\end{pf}\\

We will also use some elementary facts about sums of (imaginary parts of) logarithms, i.e.\ sums of arguments of complex numbers. Recall that we use the convention that $-\pi< arg(z)\le\pi$ for $z\in \C\setminus\left\{0\right\}$. In particular, $\log(\frac{1}{z})=-\log(z)$ for all $z\not=0$. Whenever 
$Im(z)>0$ holds, one has the equality $\log(z)=\log(-z)+\pi i$.
The following lemma collects some further elementary facts which we will use especially in the proof of \hyperref[zpq]{Lemma \ref*{zpq}}.
\begin{lem}\label{log}
For all $z\in\C\setminus\R$ we have the following identities.

$$ i)\hspace{0.3in} \log(z)-\log(\oz)=\left\{\begin{array}{cc}\log(\frac{z}{\oz})&Re(z)> 0\ \mbox{or}\ (Re(z)=0, Im(z)>0)\\
\log(\frac{z}{\oz})+2\pi i&Re(z)<0, Im(z)>0\\
\log(\frac{z}{\oz})-2\pi i&Re(z)\le 0, Im(z)<0\end{array}\right\}$$
$$ii)\hspace{0.3in}\log(z)-\log(z-\oz)=\log(\frac{z}{z-\oz}), \log(\oz)-\log(\oz-z)=\log(\frac{\oz}{\oz-z})$$
$$iii)\hspace{0.3in}\log(1-\oz)-\log(z)=\log(\frac{1-\oz}{z}), \log(1-z)-\log(\oz)=\log(\frac{1-z}{\oz})$$
$$iv)\hspace{0.3in}\log(\oz-z)-\log(\oz-1)=\log(\frac{z-\oz}{1-\oz})=\log(\oz-z)-\log(\oz-1)$$
$$v)\hspace{0.3in}\log(z)-\log(1-z)=\log(\frac{z}{1-z}), \log(\oz)-\log(1-\oz)=\log(\frac{\oz}{1-\oz})$$
$$vi)\hspace{0.3in}\log(z)+\log(1-\oz)=\log(z(1-\oz)), \log(\oz)+\log(1-z)=\log(\oz(1-z))$$
$$vii)\hspace{0.3in}\log(z(1-\oz))-\log(z-\oz)=\log(\frac{z(1-\oz)}{z-\oz}), \log(\oz(1-z))-\log(\oz-z)=\log(\frac{\oz(1-z)}{\oz-z})$$
$$viii)\hspace{0.3in}\log(1-z)-\log(1-\oz)=\left\{\begin{array}{cc}\log\frac{1-z}{1-\oz}&Re(z)<1\ \mbox{or}\ (Re(z)=1, Im(z)<0)\\
\log\frac{1-z}{1-\oz}+2\pi i&Re(z)> 1, Im(z)<0\\
\log\frac{1-z}{1-\oz}-2\pi i&Re(z)\ge 1, Im(z)>0\end{array}\right\}$$
$$ix)\hspace{0.3in}\log(\oz(1-z))-\log(z(1-\oz))=\left\{\begin{array}{cc}
\log(\frac{\oz(1-z)}{z(1-\oz)})+2\pi i&Im(z)<0\\
\log(\frac{\oz(1-z)}{z(1-\oz)})-2\pi i& Im(z)>0\end{array}\right\}$$
\end{lem}

\begin{pf} i) is obvious and ii) follows from the fact that the imaginary parts of $z$ and $z-\oz$ have the same sign, so the difference of their arguments must be smaller than $\pi$. Similarly iii) follows because the imaginary parts of $1-\oz$ and $z$ have the same sign and iv) follows because the imaginary parts of $z-\oz$ and $1-\oz$ have the same sign.

From $\mid arg(z)-arg(-z)\mid=\pi$ one can easily conclude $\mid arg(z)-arg(1-z)\mid<\pi$, which implies v). Similarly from $\mid arg(z)+arg(-\oz)\mid=\pi$ one can easily conclude $\mid arg(z)+arg(1-\oz)\mid <\pi$, which implies vi).

For vii), one can check by explicit computation that $Im(z(1-\oz))=Im(z)$, hence the imaginary parts of $z(1-\oz)$ and $z-\oz$ have the same sign and the claim follows.

viii) follows from i). For ix), assume w.l.o.g. $Im(z)>0$ and let $\alpha=arg (z)$ and $\beta=arg(1-z)$. From $arg(1-z)>arg(-z)$ we get 
$\alpha+\beta<\pi$. Then $arg(\oz(1-z))=-\alpha-\beta$ and $arg(z(1-\oz))=
\alpha+\beta$, so the difference is $-2(\alpha+\beta)$ which is smaller than $-\pi$ but bigger than $-3\pi$.
\end{pf}

\subsection{Example: the figure eight knot complement}\label{figureeight}
As the computations in the next section will be rather lengthy it shall be helpful to have an explicit example at hand to check correctness of the calculations at each step.

Consider the figure eight knot complement with its well-known ideal triangulation by two ideal simplices, let $L$ be the lift of this triangulation to the universal cover. A fundamental domain for the action of $\pi_1M$ on $L$ has 5 vertices $v_0,\ldots, v_4$, where $(v_1,v_2,v_4)$ is the common face and the (order-preserving) gluing sends $(v_0,v_1,v_2)$ to $(v_1,v_3,v_4)$, $(v_0,v_2,v_4)$ to $(v_1,v_2,v_3)$ and $(v_0,v_1,v_4)$ to $(v_2,v_3,v_4)$, see \cite[Section 4.4.2]{mr}.

With \cite[Example 3.1.1]{ggz} and the algorithm in \cite[Section 9]{gtz} we obtain 
the $PSL(2,\C)/N$-valued decoration
(with the $N$-cosets of course depending
only on the first column)
given by 

$g_{v_0} N=N$ and 
$$\hspace{-0.5in}g_{v_1} N=\left(\begin{array}{cc}0&-1\\
1&0\end{array}\right)N,g_{v_2} N=\left(\begin{array}{cc}-\omega&-\omega^2\\
-\omega&0\end{array}\right), g_{v_3} N=\left(\begin{array}{cc}-1&0\\
\omega^2-1&-1\end{array}\right)N, g_{v_4} N=\left(\begin{array}{cc}-\omega&1\\
-1&0\end{array}\right)N.$$

Let  $\Delta_{\overline{\omega}
}=(v_0,v_1,v_2,v_4)$ and $\Delta_{\omega
}=(v_1,v_2,v_3,v_4
)$, then the fundamental class for the hyperbolic monodromy $\rho$ (see \hyperref[invariant]{Definition \ref*{invariant}}) is $$\rho_*\left[M,\partial M\right]=\lambda(\Delta_{\overline{\omega}},\rho)-\lambda(\Delta_{\omega},\rho)=(\log\overline{\omega},\log\omega)-(\log\omega,\log\overline{\omega})\in\widehat{\mathcal{B}}(\C)$$ with $\omega=\frac{1}{2}+\frac{\sqrt{3}}{2}$, see \cite[Section 15]{neu}.

The Ptolemy coordinates of a $PSL(2,\C)$-representation are defined only up to sign. An obstruction cycle (for lifting 
a boundary-unipotent $PSL(2,\C)$-representation to a boundary-unipotent $SL(2,\C)$-representation, cf.\ \cite[Section 1.3]{ggz}) is given by $\sigma=(v_0,v_1,v_2)+(v_0,v_1,v_4)$. With that obstruction cycle given, we can choose the signs of the ptolemy coordinates such that the equation $\sigma_0\sigma_3c_{03}c_{12}+\sigma_0\sigma_1c_{01}c_{23}=\sigma_0\sigma_2c_{02}c_{13}$ from \cite[Definition 3.5]{ggz} is satisfied for both simplices. We obtain
so for the 
simplex $\Delta_{\overline{\omega}
}$  (with $\sigma_2=\sigma_3=-1$)
$$c_{01}=1, c_{02}=
\omega, c_{03}=1,
c_{12}=c_{13}=\omega, c_{23}=1,$$
and for the simplex 
$\Delta_{\omega
}$ (with $\sigma_0=\sigma_1=-1$)
$$c_{01}=\omega, c_{02}=1, c_{03}=\omega, c_{12}= c_{13}=1, c_{23}=\omega.$$

Now, even though this is only a $PSL(2,\C)/N$-decoration, the $\pm1$-ambiguity will disappear when we 
consider $\rho\otimes\overline{\rho}$. So (denoting by abuse of notation $N\otimes \overline{N}\subset SL(4,\C)$ again by $N$) we obtain an $SL(4,\C)/N$-decoration for $\rho\otimes\overline{\rho}$ by 

$g_{v_0} N=N$ and 
$$g_{v_1} N=\left(\begin{array}{cccc}0&0&0&1\\
0&0&-1&0\\
0&-1&0&0\\
1&0&0&0\end{array}\right)N,g_{v_2} N=\left(\begin{array}{cccc}1&\overline{\omega}&\omega&1\\
1&0&\omega&0\\
1&\overline{\omega}&0&0\\
1&0&0&0\end{array}\right)N,$$
$$ g_{v_3} N=\left(\begin{array}{cccc}1&0&0&0\\
1-\overline{\omega}^2&1&0&0\\
1-\omega^2&0&1&0\\
3&1-\omega^2&1-\overline{\omega}^2&1\end{array}\right)N, g_{v_4} N=\left(\begin{array}{cccc}1&-\omega&-\overline{\omega}&1\\
\omega&0&-1&0\\
\overline{\omega}&-1&0&0\\
1&0&0&0\end{array}\right)N.$$

One can check that for both simplices $\Delta_{\overline{\omega}}$ and $\Delta_\omega$ this is in the normal form from \hyperref[standard]{Section \ref*{standard}}, but for $\Delta_{\overline{\omega}}$ with $b$ and $c$ replaced by $-b$ and $-c$. (And for $\Delta_\omega$ one has of course to multiply by $\left(\begin{array}{cc}-1&1\\
-1&0\end{array}\right)\otimes \left(\begin{array}{cc}-1&1\\
-1&0\end{array}\right)$ from the left to bring it in normal form.) Then one may use the above Ptolemy coordinates to go through the computations in the next section and this yields (in analogy to the computations leading to the proof of \hyperref[cancel]{Corollary \ref*{cancel}}):
$$\lambda(\Delta_{\overline{\omega}},\rho\otimes\overline{\rho})=2(\log(\overline{\omega}),\log(\omega))+2(\log(\omega),\log(\overline{\omega}))$$
$$\lambda(\Delta_{\omega},\rho\otimes\overline{\rho})=2(\log(\omega),\log(\overline{\omega}))+2(\log(\overline{\omega}),\log(\omega))$$
and so $(\rho\otimes\overline{\rho})_*\left[M,\partial M\right]=\lambda(\Delta_{\overline{\omega}},\rho\otimes\overline{\rho})-\lambda(\Delta_{\omega},\rho\otimes\overline{\rho})=0$ in $\widehat{\mathcal{B}}(\C)$.

We will take up this example in \hyperref[nonlift]{Section \ref*{nonlift}}.

\section{Computations}\label{compu}
In this section we will compute the fundamental class in the extended Bloch group 
for representations of the form $\rho\otimes\overline{\rho}$ when $\rho\colon\pi_1M\to PSL(2,\C)$ is a representation 
of some $3$-manifold group. 

These representations are equivalent in $GL(4,\C)$ to those coming from the composition of $\rho$ 
with the isomorphism $PSL(2,\C)=SO(3,1)$ and it might, at first glance, have seemed more natural to 
compute the fundamental class
directly for that representation. It turns out however that that would have been much harder because a simplex 
with an (in the sense of \hyperref[ptolemy1]{Definition \ref*{ptolemy1}}) generic $PSL(2,\C)/N$-decoration need not have a 
generic $SO(3,1)/N$-decoration: in general some of the ptolemy coordinates may be zero. (One can check that the canonical triangulation of the figure eight knot complement from \hyperref[figureeight]{Section \ref*{figureeight}} yields an instance of this phenomenon.)

So further subdivision of the 
triangulation would be necessary to obtain generic $SO(3,1)/N$-decorations and this would of course vastly hamper 
computations. For this reason we will work with the representation $\rho\otimes\overline{\rho}\colon\pi_1M\to SL(4,\C)$.

\subsection{Standard form for simplices in G/N}\label{standard}
The proof of \hyperref[thm1]{Theorem \ref*{thm1}} will work by a simplexwise computation, so for most of this section we will consider one simplex and try to compute its contribution to the fundamental class and hence the Chern-Simons invariant. At first we describe a standard form for $SL(2,\C)/N$-decorated simplices in the sense of \hyperref[deco]{Definition \ref*{deco}}.

Consider $G=SL(2,\C)$ and $N\subset G$ the subset of upper triangular matrices with $1$'s on
the diagonal. 
There is a $G$-equivariant bijection $G/N=\C^2\setminus\left\{0\right\}$.
A $4$-tuple $(v_0,v_1,v_2,v_3)$ of pairwise distinct elements in $\C^2\setminus\left\{0\right\}$ is in the $SL(2,\C)$-orbit of some $4$-tuple with $v_0=(1,0)$ and $v_1=(0,a), a\in\C\setminus\left\{0\right\}$. This implies that each decorated $3$-simplex $\Delta$ is in the $SL(2,\C)$-orbit  
of a decorated simplex
$$((\begin{array}{cc}1&0\\
0&1\end{array})N
,(
\begin{array}{cc}0&-\frac{1}{a}\\
a&0\end{array})N
,(\begin{array}{cc}-\frac{d}{a}&-\frac{1}{b}\\
b&0\end{array})N
,(\begin{array}{cc}-\frac{e}{a}&-\frac{1}{c}\\
c&0\end{array})N).
$$
\begin{tikzpicture}
\path[draw](3,3)--(7.5,1.2)--(6,0);
\path[draw, fill=green!20](0,0)--(6,0)--(7.5,1.2);
\path[draw](0,0)--(3,3)--(6,0);
\draw[dashed](0,0)--(7.5,1.2);
\fill (1.5,1.5) circle (1.6pt) node [below=5pt] {a};
\fill (3,0) circle (1.6pt) node [below=5pt] {b};\fill (3.75,0.6) circle (1.6pt) node [below=5pt] {c};
\fill (4.5,1.5) circle (1.6pt) node [below=5pt] {d};\fill (5.25,2.1) circle (1.6pt) node [below=5pt] {e};
\fill (6.75,0.6) circle (1.6pt) node [below=5pt] {f};
\fill (0,0) circle (1.6pt) node [below=5pt] {$\left(\begin{array}{cc}1&0\\
0&1\end{array}\right)N$};
\fill (3,3) circle (1.6pt) node [above=3pt] {$\left(\begin{array}{cc}0&-\frac{1}{a}\\
a&0\end{array}\right)N$};
\fill (6,0) circle (1.6pt) node [below=5pt] {$\left(\begin{array}{cc}-\frac{d}{a}&-\frac{1}{b}\\
b&0\end{array}\right)N$};\fill (7.5,1.2) circle (1.6pt) node [right=4pt] {$\left(\begin{array}{cc}-\frac{e}{a}&-\frac{1}{c}\\
c&0\end{array}\right)N$};
\end{tikzpicture}

One may check that the Ptolemy coordinates (in the sense of \hyperref[ptolemy1]{Definition \ref*{ptolemy1}}) of the edges are 
$$c_{1100}=a, c_{1010}=b, c_{1001}=c, c_{0110}=d, c_{0101}=e,c_{0011}=f$$
with $af+cd=be$,
and thus $$\lambda(c_{0000})
=(\log(c)+\log(d)-\log(b)-\log(e),\log(a)+\log(f)-\log(b)-\log(e)).$$

When $\Delta\subset H^3\cup\partial_\infty H^3$ is the ideal hyperbolic simplex whose ideal vertices are the projections of $v_0,v_1,v_2,v_3$ to $\C P^1=\partial_\infty H^3$, then one can easily check that $z=\frac{cd}{be}$ is the cross ratio of the vertices of $\Delta$ and it is well-known that $\Delta$ is non-degenerate if and only if $z\not\in\R$ and that the ordering of $\Delta$ agrees with the orientation of $H^3$ if and only if $Im(z)>0$.

We remark for later use that 
$\lambda(c_{0000})=(\log(z)+2p\pi i,\log(1-z)+2q\pi i)$ for some integers $p$ and $q$. 

We will see in the next subsection that our computations will only work for $z\not\in\R$, i.e., for non-degenerate simplices (although in the definition of the Bloch groups only $z\not\in\left\{0,1\right\}$ is required). At the time of writing it is not known whether every hyperbolic $3$-manifold admits an ideal triangulation with no degenerate simplex. However the methods of \cite{gtz} do not require ideal triangulations but allow interior vertices, so upon performing barycentric subdivision and suitably decorating the interior vertices we can always assume to have simplices with cross ratios $z\not\in \R$ throughout. (See \cite[Proposition 5.4]{gtz}.)


\subsection{Toy case: computations in the (non-extended) Bloch group.} Because the following computations in the extended (pre-)Bloch group might seem a bit unmotivated at first glance, we start this section with explaining the proof of a considerable simpler fact, namely the equality $$(\rho\otimes\overline{\rho})_*\left[M,\partial M\right]=2\rho_*\left[M,\partial M\right]+
2\overline{\rho_*\left[M,\partial M\right]}\in{\mathcal{B}}(\C)$$ in the (non-extended) Bloch group ${\mathcal{B}}(\C)$. (In terms of the Cheeger-Chern-Simons invariant this means the vanishing of its imaginary part $Vol(\rho)$, which of course is already implied by the complex conjugacy invariance from  \hyperref[conj]{Lemma \ref*{conj}}.) We hope that these computations in the (pre-)Bloch group
prepare the reader for the more complicated computations in the extended (pre-)Bloch group which will be performed in the remainder of this section.

So let $\Delta$ be an ideal simplex in a triangulation of a $3$-manifold. (Let us stick to ideal triangulations for simplicity of exposition, although this restriction is not necessary for the argument.) The boundary-parabolic representation $\rho\colon\pi_1M\to SL(2,\C)$ comes with an equivariant map $\widetilde{M}\to H^3$ which maps a lift of $\Delta$ to a hyperbolic ideal simplex.
Its cross ratio defines an element $\left[z\right]\in{\mathcal P}(\C)$ in the pre-Bloch group 
and the fundamental class $\rho_*\left[M,\partial M\right]$ is the sum of these elements over all simplices in the triangulation.

Similarly one can compute the contribution of $\Delta$ to the 
fundamental class $\rho\otimes\overline{\rho}\left[M,\partial M\right]$ from $\left[z\right]\in{\mathcal P}(\C)$ 
alone. Namely, \hyperref[gtz]{Proposition \ref*{gtz}} tells us that we have to formally consider 10 subsimplices of $\Delta$ 
and look at their contribution. The computations for this are in principle the same as 
in the following \hyperref[tensor]{Section \ref*{tensor}}, but they are simpler 
because 
we do not need to know the elements 
$\lambda(c_{ijkl})\in\PC(\C)$, but rather only their associated cross ratios 
$cr(\lambda(c_{ijkl}))\in{\mathcal P}(\C)$. Those are 
with $z=\frac{cd}{be}=1-\frac{af}{be}$ easily computed as follows:
$$cr(\lambda(c_{2000}))=cr(\lambda(c_{0200}))=cr(\lambda(c_{0020}))=cr(\lambda(c_{0002}))=\left[\oz\right]$$
$$cr(\lambda(c_{1100}))=cr(\lambda(c_{0011}))=\left[\frac{z}{\oz}\right]$$
$$cr(\lambda(c_{1010}))=cr(\lambda(c_{0101}))=\left[\frac{z(1-\oz)}{z-\oz}\right]$$
$$cr(\lambda(c_{1001}))=cr(\lambda(c_{0110}))=\left[\frac{z-\oz}{1-\oz}\right],$$
and so we obtain the contribution of this simplex to the fundamental class of $\rho\otimes\overline{\rho}$ in ${\mathcal{B}}(\C)$ to be
$$4\left[\oz\right]+2 \left[\frac{z}{\oz}\right]+2\left[\frac{z(1-\oz)}{z-\oz}\right]+2\left[\frac{z-\oz}{1-\oz}\right].$$

Now in ${\mathcal P}(\C)$ a direct application of the 5-term relation yields an equality
$$\left[\frac{z-\oz}{1-\oz}\right]-\left[z\right]+\left[\frac{z(1-\oz)}{z-\oz}\right]-\left[\frac{z}{z-\oz}\right]+
\left[\frac{1}{1-\oz}\right]=0,$$
which we can simplify by using the 
following relations, which actually are consequences of well-known functional equations for the Bloch-Wigner dilogarithm $D_2$:
$$\left[\frac{1}{1-\oz}\right]=\left[\oz\right], \left[\frac{z}{z-\oz}\right]=\left[\frac{1}{1-\frac{\oz}{z}}\right]=\left[\frac{\oz}{z}\right]=
-\left[\frac{z}{\oz}\right],
$$ 
so that we obtain  
$$\left[\frac{z-\oz}{1-\oz}\right]-\left[z\right]+\left[\frac{z(1-\oz)}{z-\oz}\right]+\left[\frac{z}{\oz}\right]+\left[\oz\right]=0.$$

Plugging twice of this in the above formula we obtain that the contribution of the simplex to the fundamental class of $\rho\otimes\overline{\rho}$ simplifies to
$$2\left[\oz\right]+2\left[z\right],$$
and this implies the claimed equality in ${\mathcal{B}}(\C)$.

To carry through an analogous computation in the extended pre-Bloch group we will have to replace the cross ratios by the elements $\lambda(c_{ijkl})\in\PC(\C)$. This will lead to considerably more complicated computations, but we will in an analogous way as above use the extended 5-term relation and the symmetries of the extended Rogers' dilogarithm which we derived in \hyperref[gz]{Lemma \ref*{gz}}, to simplify the formulas and will so in the end obtain the equality of \hyperref[thm1]{Theorem \ref*{thm1}} also in $\widehat{\mathcal{B}}(\C)$.

\subsection{Ptolemy coordinates of the tensor product}\label{tensor}
For a representation $$\rho
\colon\pi_1M\to SL(2,\C)$$
we consider its tensor product with its complex conjugate
$$\rho\otimes\overline{\rho}\colon\pi_1M\to SL(4,\C).$$ Triangular matrices with 1's on the diagonal are sent to triangular matrices with 1's on the diagonal.

Given an $SL(2,\C)/N$-decoration for $\rho$, we obtain a decoration for $\rho\otimes\overline{\rho}$ in the obvious way, replacing $A$ by $A\otimes\overline{A}$ on each vertex. 

So we fix again one simplex $\Delta$, then the normal form decoration from the previous \hyperref[standard]{Section \ref*{standard}} is mapped to the $N$-cosets of
$$\hspace{-0.5in}\left(\left(
\begin{array}{cccc}
1&0&0&0\\
0&1&0&0\\
0&0&1&0\\
0&0&0&1
\end{array}\right)
,\left(
\begin{array}{cccc}
0&0&0&\frac{1}{\mid a\mid^2}\\
0&0&-\frac{\oa}{a}&0\\
0&-\frac{a}{\oa}&0&0\\
\mid a\mid^2&0&0&0\end{array}\right)
,\left(
\begin{array}{cccc}
\mid\frac{d}{a}\mid^2&\frac{d}{a\ob}&\frac{\od}{b\oa}&\frac{1}{\mid b\mid^2}\\
-\frac{d\ob}{a}&0&-\frac{\ob}{b}&0\\
-\frac{b\od}{\oa}&-\frac{b}{\ob}&0&0\\
\mid b\mid^2&0&0&0\end{array}\right)
,\left(
\begin{array}{cccc}
\mid\frac{e}{a}\mid^2&\frac{e}{a\oc}&\frac{\oee}{c\oa}&\frac{1}{\mid c\mid^2}\\
-\frac{e\oc}{a}&0&-\frac{\oc}{c}&0\\
-\frac{c\oee}{\oa}&-\frac{c}{\oc}&0&0\\
\mid c\mid^2&0&0&0\end{array}\right)\right)
$$
From this and $af+cd=be$ we compute the ptolemy coordinates of $\Delta$ as follows
$$
c_{3100}=\mid a\mid^2,
c_{3010}=\mid b\mid^2,
c_{3001}=\mid c\mid^2$$
$$
c_{2200}=a^2,
c_{2110}=ab\od,
c_{2101}=ac\oee$$
$$c_{2020}=b^2,
c_{2011}=bc\of,
c_{2002}=c^2$$
$$
c_{1300}=-\mid a\mid^2,
c_{1210}=-a\ob d,
c_{1201}=-a\oc e$$
$$
c_{1120}=\oa bd,
c_{1111}=2Im(\ob cd\oee)i,
c_{1102}=\oa ce$$
$$
c_{1030}=-\mid b\mid^2,
c_{1021}=-b\oc f,
c_{1012}=\ob cf$$
$$
c_{1003}=-\mid c\mid^2,
c_{0310}=\mid d\mid^2,
c_{0301}=\mid e\mid^2$$
$$
c_{0220}=d^2,
c_{0211}=de\of,
c_{0202}=e^2$$
$$
c_{0130}=-\mid d\mid^2,
c_{0121}=-d\oee f                                                                                                                                                                                                                                                                                                                                                                                                                                                                                                                                                                          ,
c_{0112}=\od ef$$
$$
c_{0103}=-\mid e\mid^2,
c_{0031}=\mid f\mid^2,
c_{0022}=f^2,
c_{0013}=-\mid f\mid^2
$$
(It should have sufficed to compute $c_{3100},c_{2200},c_{2110},c_{1111}$ and then proceed by symmetry, however to be safe we doublechecked all computations with \cite{sage}.)

We note that these coordinates are nonzero (i.e., the decoration is generic in the sense of 
\hyperref[generi]{Definition \ref*{generi}}) if and only if the Ptolemy coordinates $a,\ldots,f$ of $\rho$ are nonzero and if moreover $\ob cd\oee \not\in\R$. The latter condition is in view of $\ob cd\oee=\mid be\mid^2\frac{cd}{be}$ equivalent to the condition that the cross ratio $z=\frac{cd}{be}$ is not a real number and as argued in \hyperref[standard]{Section \ref*{standard}} above this can always be assumed.

We plug the ptolemy coordinates of $\Delta$
into the formula from \hyperref[invariant]{Definition \ref*{invariant}} and obtain the following.
$$\hspace{-0.5in}\lambda(c_{2000})=(\log\mid c\mid^2+\log ab\od-\log\mid b\mid^2-\log ac\oee,\log\mid a\mid^2+\log bc\of-\log\mid b\mid^2-\log ac\oee)$$
$$\hspace{-0.5in}\lambda(c_{1100})=(\log ac\oee+\log(-a\ob d)-\log ab\od-\log(-a\oc e),\log a^2+\log 2Im(\ob cd\oee)i-\log ab\od-\log(-a\oc e)$$
$$\hspace{-0.5in}\lambda(c_{1010})=(\log bc\of+\log\oa bd-\log b^2-\log 2Im(\ob cd\oee)i,\log ab\od+\log(-b\oc f)-\log b^2-\log 2Im(\ob cd\oee)i$$
$$\hspace{-0.5in}\lambda(c_{1001})=(\log c^2+\log 2Im(\ob cd\oee)i-\log bc\of-\log \oa ce,\log ac\oee+\log \ob cf-\log bc\of-\log \oa ce)$$
$$\hspace{-0.5in}\lambda(c_{0200})=(\log(-a\oc e)+\log\mid d\mid^2-\log(-a\ob d)-\log\mid e\mid^2,\log(-\mid a\mid^2)+\log de\of-\log(-a\ob d)-\log\mid e\mid^2)$$
$$\hspace{-0.5in}\lambda(c_{0110})=(\log 2Im(\ob cd\oee)i+\log d^2-\log\oa bd-\log de\of,\log(-a\ob d)+\log (-d \oee f)-\log\oa bd-\log de\of)$$
$$\hspace{-0.5in}\lambda(c_{0101})=(\log \oa ce+\log de\of-\log 2Im(\ob cd\oee)i-\log e^2,\log(-a\oc e)+\log \od ef-\log 2Im(\ob cd\oee)i-\log e^2)$$
$$\hspace{-0.5in}\lambda(c_{0020})=(\log (-b\oc f)+\log(-\mid d\mid^2)-\log (-\mid b\mid^2)-\log(-d\oee f),\log \oa bd+\log \mid f\mid^2-\log (-\mid b\mid^2)-\log(-d\oee f))$$
$$\hspace{-0.5in}\lambda(c_{0011})=(\log \ob cf+\log (-d\oee f)-\log(-b\oc f)-\log \od ef,\log 2Im(\ob cd\oee)i+\log f^2-\log(-b\oc f)-\log \od ef)$$
$$\hspace{-0.5in}\lambda(c_{0002})=(\log(-\mid c\mid^2)+\log\od ef-\log\ob cf-\log (-\mid e\mid^2),\log\oa ce+\log (-\mid f\mid^2)-\log\ob cf-\log (-\mid e\mid^2))$$

By \hyperref[lift]{Proposition \ref*{lift}} we are free to chose lifts (i.e., logarithms) as long as we chose the same lift for the same numbers. To simplify the above expressions we fix once and for all the following lifts:
$$\log(ab\od)=\log a+\log b+\log\od, \log a\ob d=\log a+\log\ob +\log d, etc.$$
$$\log(-ab\od)=\pi i+\log a+\log b+\log\od, etc.$$
$$\log a^2=2\log a, log(-a^2)=\pi i+2\log a, etc.$$
$$\log\mid a\mid^2=\log a+\log\oa, \log(-\mid a\mid^2)=\pi i+\log a+\log\oa, etc.$$
$$\log(2Im(\ob cd\oee)i)=\pi i+\log b+\log\oc+\log\od+\log e+\log(1-\frac{\ob cd\oee}{b\oc\od e})$$
The reason for the somewhat unnatural seeming choice of $\tilde{\lambda}(c_{1111})$ is to simplify the formulas in the following lemma. One can easily check that it is indeed a choice of logarithm, i.e., that exponentiating the right hand side gives $2Im(\ob cd\oee)i$.

Now with these lifts we can simplify as follows. 
\begin{lem}\label{ptolemy}With the above-chosen lifts, the pre-Bloch group elements associated to the Ptolemy coordinates of $\rho\otimes\overline{\rho}$ for the decorated simplex $\Delta$ from \hyperref[standard]{Section \ref*{standard}} are:
$$\hspace{-3in}\tilde{\lambda}(c_{2000})=\tilde{\lambda}(c_{0200})=\tilde{\lambda}(c_{0002})=$$
$$\hspace{1in}(-\log\ob+\log\oc+\log\od-\log\oee,\log\oa-\log\ob-\log\oee+\log\of)$$
$$\tilde{\lambda}(c_{0020})=(-\log\ob+\log\oc+\log\od-\log\oee,-2\pi i+\log\oa-\log\ob -\log\oee+\log\of)$$
$$\hspace{-3.7in}\tilde{\lambda}(c_{1100})=\tilde{\lambda}(c_{0011})=$$
$$\hspace{0.5in}(\log\ob-\log b+\log c-\log\oc +\log d-\log\od+\log\oee-\log e,\log(1-\frac{\ob cd\oee}{b\oc\od e}))$$
$$\begin{array}{c}\tilde{\lambda}(c_{1010})=\tilde{\lambda}(c_{0101})=
(-\pi i+\log\oa-\log b+\log c-\log\oc+\log d-\log\od-\log e+\log\of-\log(1-\frac{\ob cd\oee}{b\oc\od e}),\\
\hspace{1in}\log a-\log b-\log e+\log f-\log(1-\frac{\ob cd\oee}{b\oc\od e}))\end{array}$$
$$\begin{array}{c}\hspace{-1in}\tilde{\lambda}(c_{1001})=(\pi i-\log\oa+\log \oc+\log \od-\log\of+\log(1-\frac{\ob cd\oee}{b\oc\od e}),\\
\hspace{1in}\log a-\log\oa+\log\ob-\log b +\log\oee-\log e+\log f-\log\of)\end{array}$$
$$\begin{array}{c}\hspace{-1in}\tilde{\lambda}(c_{0110})=(\pi i-\log\oa+\log \oc+\log \od-\log\of+\log(1-\frac{\ob cd\oee}{b\oc\od e}),\\
\hspace{1in}2\pi i+\log a-\log\oa+\log\ob-\log b +\log\oee-\log e+\log f-\log\of)\end{array}$$
\end{lem}

Remark: The symmetry-breaking formulas for $\tilde{\lambda}(c_{0020})$ and $\tilde{\lambda}(c_{0110})$ seem to be an artefact of our somewhat arbitrary choice of logarithms.

Next we want to express these formulas in the $(z,p,q)$-form (in the notation of \cite{gtz}, compare the remark after \hyperref[extprebloch]{Definition \ref*{extprebloch}}). The proof will use some elementary facts about complex logarithms, which for better readability had been collected in \hyperref[log]{Lemma \ref*{log}} before. 

\begin{lem}\label{zpq}Let $z=\frac{cd}{be}=1-\frac{af}{be}\not\in\R$ be the cross ratio of $\Delta$, and 
define the integers $p,q$ via 
$$\log c+\log d-\log b-\log e =\log(z)+2p\pi i,$$
$$\log a+\log f-\log b-\log e=\log(1-z)+2q\pi i.$$
Then
$$\tilde{\lambda}(c_{2000})=\tilde{\lambda}(c_{0200})=\tilde{\lambda}(c_{0002})=(\log(\oz)-2p\pi i,\log(1-\oz)-2q\pi i)$$
$$\tilde{\lambda}(c_{0020})=(\log(\oz)-2p\pi i,\log(1-\oz)-2(q+1)\pi i)$$
$$\tilde{\lambda}(c_{1100})=\tilde{\lambda}(c_{0011})=\left\{\begin{array}{cc}(\log(\frac{z}{\oz})+4p\pi i,\log(1-\frac{z}{\oz}))&Re(z)> 0\ \mbox{or}\ (Re(z)=0, Im(z)>0)\\
(\log(\frac{z}{\oz})+(4p+2)\pi i,\log(1-\frac{z}{\oz}))&Re(z)<0, Im(z)>0\\
(\log(\frac{z}{\oz})+(4p-2)\pi i,\log(1-\frac{z}{\oz}))&Re(z)\le 0, Im(z)<0\end{array}\right\}$$
$$\tilde{\lambda}(c_{1010})=\tilde{\lambda}(c_{0101})=\left\{\begin{array}{cc}
(\log(\frac{z(1-\oz)}{z-\oz})+2(2p-q)\pi i,\log(\frac{(1-z)\oz}{\oz-z})+2q\pi i)
&Im (z)>0\\
(\log(\frac{z(1-\oz)}{z-\oz})+2(2p-q-1)\pi i,\log(\frac{(1-z)\oz}{\oz-z})+2q\pi i)&Im(z)<0
\end{array}\right\}$$
$$\tilde{\lambda}(c_{1001})=\left\{\begin{array}{cc}(\log(\frac{\oz-z}{\oz-1})+2(q-p)\pi i,\log\frac{1-z}{1-\oz}+4q\pi i)&Im(z)>0,Re(z)<1\\
\log(\frac{\oz-z}{\oz-1})+2(q-p)\pi i,\log\frac{1-z}{1-\oz}+(4q+2)\pi i)&Im(z)>0,Re(z)\ge1\\
\log(\frac{\oz-z}{\oz-1})+2(q-p+1)\pi i,\log\frac{1-z}{1-\oz}+4q\pi i)&Im(z)<0,Re(z)\le 1\\
\log(\frac{\oz-z}{\oz-1})+2(q-p+1)\pi i,\log\frac{1-z}{1-\oz}+(4q+2)\pi i)&Im(z)<0,Re(z)> 1\end{array}\right\}$$
$$\tilde{\lambda}(c_{0110})=\left\{\begin{array}{cc}(\log(\frac{\oz-z}{\oz-1})+2(q-p)\pi i,\log\frac{1-z}{1-\oz}+(4q+2)\pi i)&Im(z)>0,Re(z)<1\\
\log(\frac{\oz-z}{\oz-1})+2(q-p)\pi i,\log\frac{1-z}{1-\oz}+(4q+4)\pi i)&Im(z)>0,Re(z)\ge1\\
\log(\frac{\oz-z}{\oz-1})+2(q-p+1)\pi i,\log\frac{1-z}{1-\oz}+(4q+2)\pi i)&Im(z)<0,Re(z)\le 1\\
\log(\frac{\oz-z}{\oz-1})+2(q-p+1)\pi i,\log\frac{1-z}{1-\oz}+(4q+4)\pi i)&Im(z)<0,Re(z)> 1\end{array}\right\}$$
\end{lem}

\begin{pf}
a) The formulas for $\tilde{\lambda}(c_{2000})=\tilde{\lambda}(c_{0200})=\tilde{\lambda}(c_{0002})$ and $\tilde{\lambda}(c_{2000})$ are immediate by complex conjugation of the formulas for $\log(z)+2p\pi i$ and $\log(1-z)+2q\pi i$.

b) Subtraction yields $\tilde{\lambda}(c_{1100})=\tilde{\lambda}(c_{0011})=(\log(z)-\log(\oz)+4p\pi i,\log(1-\frac{z}{\oz}))$. The formulas for $\tilde{\lambda}(c_{1100})=\tilde{\lambda}(c_{0011})$ then follow from \hyperref[log]{Lemma \ref*{log}} i).

c) By \hyperref[log]{Lemma \ref*{log}} v) we have $\log(z)-\log(1-z)=\log(\frac{z}{1-z})$. So subtraction yields 
$$-\log(a)+\log(c)+\log(d)-\log(f)=\log(\frac{z}{1-z})+2(p-q)\pi i.$$
From this one obtains
$$\begin{array}{c}\hspace{-3in}\tilde{\lambda}(c_{1010})=\tilde{\lambda}(c_{0101})=\\
(-\pi i-\log(\frac{\oz}{1-\oz})+2(p-q)\pi i +\log(z)+2p\pi i-\log(1-\frac{z}{\oz}),\log(1-z)+2q\pi i-\log(1-\frac{z}{\oz}))\end{array}$$
Using ii) and v) from \hyperref[log]{Lemma \ref*{log}} the first coordinate simplifies to $-\pi i+\log(1-\oz)+\log(z)-\log(\oz-z)+2(2p-q)\pi i$. We observe that $-\log(\oz-z)-\pi i$ is $-\log(z-\oz)$ or $-\log(z-\oz)-2\pi i$ according to whether $Im(z)>0$ or $Im(z)<0$. Then with \hyperref[log]{Lemma \ref*{log}} vi) and vii) we get the claimed formula for the first coordinate. The formula for the second coordinate is also obtained by applying \hyperref[log]{Lemma \ref*{log}} vi) and vii).

d) Conjugating the first equation in the proof of c) and plugging it into the first coordinate yields $$\tilde{\lambda}(c_{1001})=(\pi i+\log(\frac{\oz}{1-\oz})-2(p-q)\pi i+\log(1-\frac{z}{\oz}),\log(1-z)+2q\pi i-\log(1-\oz)+2q\pi i)$$
We have again $\log(\frac{\oz}{1-\oz})=\log(\oz)-\log(1-\oz)$ and $\log(1-\frac{z}{\oz})=\log(\oz-z)-\log(\oz)$.  If $Im(z)>0$, 
then $\pi i+\log(\oz)=\log(-\oz)$ and (as a direct computation 
shows) $Im(\frac{\oz}{1-\oz})<0$, so $\pi i +
\log(\frac{\oz}{1-\oz})=\log(\frac{-\oz}{1-\oz})$, so the first coordinate simplifies to $\pi i-\log(1-\oz)+\log(\oz-z)+2(q-p)\pi i=-\log(\oz-1)+\log(\oz-z)+2(q-p)\pi i=\log(\frac{\oz-z}{\oz-1})+2(q-p)\pi i$,
where the last equality uses that $z-\oz$ and $1-\oz$ have the same imaginary part, 
so the difference of their arguments is in $(-\pi,\pi)$. If $Im(z)<0$, then one obtains by similar arguments that the first coordinate equals $\log(\frac{\oz-z}{\oz-1})+2(q-p+1)\pi i$.
The formula for the second coordinate follows from \hyperref[log]{Lemma \ref*{log}} viii).
This proves the formula for $\tilde{\lambda}(c_{1001})$ and by exactly the same arguments we obtain that for $\tilde{\lambda}(c_{0110})$.
\end{pf}

\begin{cor}\label{corsum}Under the assumptions of \hyperref[zpq]{Lemma \ref*{zpq}} we have
$$\tilde{\lambda}(c_{2000})+\tilde{\lambda}(c_{0200})+\tilde{\lambda}(c_{0002})+\tilde{\lambda}(c_{0020})+\tilde{\lambda}
(c_{1100})+\tilde{\lambda}(c_{0011})+\tilde{\lambda}(c_{1010})+\tilde{\lambda}(c_{0101})+\tilde{\lambda}(c_{1001})
+\tilde{\lambda}(c_{0110})=$$
$$4(\log(\oz)-2p\pi i,\log(1-\oz)-2q\pi i)$$
$$+2(\log(\frac{z}{\oz})+4p\pi i,\log(1-\frac{z}{\oz}))$$
$$+2(\log(\frac{z(1-\oz)}{z-\oz})+2(2p-q)\pi i,\log(\frac{(1-z)\oz}{\oz-z})+2q\pi i)$$
$$+2(\log(\frac{\oz-z}{\oz-1})+2(q-p)\pi i,
\log(\frac{1-z}{1-\oz})+2(2q+1)\pi i)$$
$$+\left\{\begin{array}{ccc}
\chi(-\log(\oz)+\log(\frac{\oz-z}{\oz-1})+2q\pi i)
&:&Re(z)\ge 1, Im(z)>0\\
\chi(-\log(\oz)-\log(\frac{\oz-z}{\oz-1})+2q\pi i)
&:&1>Re(z)\ge 0, Im(z)>0\\
\chi(-\log(\oz)+2\log(1-\frac{z}{\oz})-\log(\frac{\oz-z}{\oz-1})+2q\pi i)
&:&Re(z)< 0, Im(z)>0\\
\chi(-\log(\oz)+2\log(\frac{(1-z)\oz}{\oz-z})-2\log(\frac{1-z}{1-\oz})+\log(\frac{\oz-z}{\oz-1})+2(q-1)\pi i)
&:&Re(z)> 1 , Im(z)<0\\
\chi(-\log(\oz)+2\log(\frac{(1-z)\oz}{\oz-z})-2\log(\frac{1-z}{1-\oz})-\log(\frac{\oz-z}{\oz-1})+2(q-1)\pi i)
&:&1\ge Re(z)> 0, Im(z)<0\\
\chi(-\log(\oz)-2\log(1-\frac{z}{\oz})+2\log(\frac{(1-z)\oz}{\oz-z})-2\log(\frac{1-z}{1-\oz})-\log(\frac{\oz-z}{\oz-1})+2(q-1)\pi i)
&:&Re(z)\le 0, Im(z)<0
\end{array}\right\}$$
\end{cor}
\begin{pf}This follows from \hyperref[zpq]{Lemma \ref*{zpq}} by using part iv) of \hyperref[hom]{Lemma \ref*{hom}} and $\chi(4\pi i)=0$.\end{pf}

\subsection{Using the five-term relation}
\begin{lem}\label{blochsum}In $\PC(\C)$ we have the equality 
$$\left\{\begin{array}{ccc}
(\log(\frac{z-\oz}{1-\oz})+2(q-p)\pi i, \log(\frac{1-z}{1-\oz})+4q\pi i)
&:&Re(z)<1\ \mbox{or}\ (Re(z)=1, Im(z)<0)\\
(\log(\frac{z-\oz}{1-\oz})+2(q-p)\pi i, \log(\frac{1-z}{1-\oz})+2(2q+1)\pi i)&:&Re(z)> 1, Im(z)<0\\
(\log(\frac{z-\oz}{1-\oz})+2(q-p)\pi i, \log(\frac{1-z}{1-\oz})+2(2q-1)\pi i)&:&Re(z)\ge 1, Im(z)>0\end{array}\right\}$$
$$-(\log(z)+2p\pi i,\log(1-z)+2q\pi i)$$
$$+(\log(\frac{z(1-\oz)}{z-\oz})+2(2p-q)\pi i,\log(\frac{\oz(1-z)}{\oz-z})+2q\pi i)$$
$$-(\log(\frac{z}{z-\oz})+4p\pi i,\log(\frac{\oz}{\oz-z}))$$
$$+(\log(\frac{1}{1-\oz})+2q\pi i,\log(\frac{\oz}{\oz-1})+2(q-p)\pi i)$$
$$=0$$

for any $p,q\in \Z$.\end{lem}
\begin{pf}We want to apply the five-term relation from \hyperref[extprebloch]{Definition \ref*{extprebloch}} 
with the five summands above corresponding to $(e_0,f_0),\ldots,(e_4,f_4)$. According 
to \hyperref[extprebloch]{Definition \ref*{extprebloch}} we have to check the conditions 
$$e_2=e_1-e_0,e_3=e_1-e_0-f_1+f_0,f_3=f_2-f_1,e_4=f_0-f_1,f_4=f_2-f_1+e_0.$$
In each of the five cases the equality is a direct consequence of \hyperref[log]{Lemma \ref*{log}}.

\end{pf}

\begin{cor}\label{corfive}
In $\PC(\C)$ we have the equality 
$$
(\log(\frac{z-\oz}{1-\oz})+2(q-p)\pi i, \log(\frac{1-z}{1-\oz})+2(2q+1)\pi i)
$$

$$+(\log(\frac{z(1-\oz)}{z-\oz})+2(2p-q)\pi i,\log(1-\oz)+2q\pi i)$$
$$+(\log(\frac{z}{\oz})+4p\pi i,\log(\frac{\oz-z}{\oz}))=$$
$$(\log(z)+2p\pi i,\log(1-z)+2q\pi i)$$
$$-(\log(\oz)-2p\pi i,\log(1-\oz)-2q\pi i)$$
$$\left\{\begin{array}{ccc}0&:&Re(z)<1\ \mbox{or}\ (Re(z)=1, Im(z)<0)\\
\chi(\log(\frac{z-\oz}{1-\oz})+2(q-p)\pi i)&:&Re(z) > 1, Im(z)<0\\
-\chi(\log(\frac{z-\oz}{1-\oz})+2(q-p)\pi i)&:&Re(z)\ge 1, Im(z)>0\end{array}\right\}+$$
$$\left\{\begin{array}{ccc}
\chi(-\frac{1}{2}\log(1-\oz)-\frac{1}{2}\log(1-\frac{\oz}{z})-\frac{1}{2}\log(\frac{z}{\oz})-(2q^2+q)\pi i)&:&Re(z)>0,Im(z)>0\\
\chi(-\frac{1}{2}\log(1-\oz)+\frac{1}{2}\log(1-\frac{\oz}{z})+\frac{1}{2}\log(\frac{z}{\oz})-(2q^2+q)\pi i))&:&Re(z)< 0, Im(z)>0\\
\chi(\frac{1}{2}\log(1-\oz)+\frac{1}{2}\log(1-\frac{\oz}{z})+\frac{1}{2}\log(\frac{z}{\oz})+(2q^2+q)\pi i))&:&Re(z)> 0,Im(z)<0\\
\chi(\frac{1}{2}\log(1-\oz)-\frac{1}{2}\log(1-\frac{\oz}{z})-\frac{1}{2}\log(\frac{z}{\oz})+(2q^2+q)\pi i))&:&Re(z)<0,Im(z)<0
\end{array}\right\}$$

for any $p,q\in \Z$.\end{cor}
\begin{pf}We will prove this by adapting some of the terms in \hyperref[blochsum]{Lemma \ref*{blochsum}}. 

The correction term for the first summand is a direct application of \hyperref[hom]{Lemma \ref*{hom}}. 
Let us look at the fifth and fourth summand from \hyperref[blochsum]{Lemma \ref*{blochsum}}.

From \hyperref[gz]{Lemma \ref*{gz}} or \hyperref[coro]{Corollary \ref*{coro}}, respectively, we obtain
$$(\log(\frac{1}{1-\oz})+2q\pi i,\log(\frac{\oz}{\oz-1})+2(q-p)\pi i)=$$
$$\left\{\begin{array}{ccc}(\log(\oz)-2p\pi i,\log(1-\oz)-2q\pi i)+
\chi((-\frac{1}{2}\log(1-\oz)-(2q^2+q-\frac{1}{6})\pi i&:&Im(z)>0\\
(\log(\oz)-2p\pi i,\log(1-\oz)-2q\pi i)+\chi((\frac{1}{2}\log(1-\oz)+(2q^2+q+\frac{1}{6})\pi i)&:&Im(z)<0\end{array}\right\}$$

and, using $\frac{z}{z-\oz}=\frac{1}{1-\frac{\oz}{z}}$ (and an explicit computation if $Re(z)=0$, i.e., if $\frac{z}{\oz}=-1$) we get
$$(\log(\frac{z}{z-\oz})+4p\pi i,\log(\frac{\oz}{\oz-z}))=$$
$$(\log(\frac{\oz}{z})-4p\pi i,\log(1-\frac{\oz}{z})-4p\pi i)+\left\{\begin{array}{ccc}
\chi(\frac{1}{2}\log(1-\frac{\oz}{z})+(2p+\frac{1}{6})\pi i)&:&Re(z)\ge 0,Im(z)>0\\
&&
\mbox{ or }Re(z)\le 0,Im(z)<0\\
\chi(-\frac{1}{2}\log(1-\frac{\oz}{z})-(2p-\frac{1}{6})\pi i)&:&Re(z)> 0,Im(z)<0\\
&&\mbox{ or }Re(z)<0, Im(z)>0\end{array}\right\}$$
$$=-(\log(\frac{z}{\oz})+4p\pi i,\log(1-\frac{z}{\oz}))+\left\{\begin{array}{ccc}
\chi(-\frac{1}{2}\log(1-\frac{\oz}{z})-\frac{1}{2}\log(\frac{z}{\oz})+\frac{\pi i}{6}))&:&Re(z)\ge 0,Im(z)>0\\
&&\mbox{ or }Re(z)\le 0,Im(z)<0\\
\chi(\frac{1}{2}\log(1-\frac{\oz}{z})+\frac{1}{2}\log(\frac{z}{\oz})+\frac{\pi i}{6}))&:&Re(z)> 0,Im(z)<0\\
&&\mbox{ or }
Re(z)< 0, Im(z)>0
\end{array}\right\}$$
Plugging this into the equation from \hyperref[blochsum]{Lemma \ref*{blochsum}} we get the claim, using
that $\chi$ is a homomorphism vanishing on multiples of $4\pi i$.
\end{pf}



\begin{cor}\label{cancel}With the notation from \hyperref[zpq]{Lemma \ref*{zpq}} we have $$\tilde{\lambda}(c_{2000})+\tilde{\lambda}(c_{0200})+\tilde{\lambda}(c_{0002})+
\tilde{\lambda}(c_{0020})+\tilde{\lambda}(c_{1100})+$$
$$+\tilde{\lambda}(c_{0011})+
\tilde{\lambda}(c_{1010})+\tilde{\lambda}(c_{0101})+
\tilde{\lambda}(c_{1001})+
\tilde{\lambda}(c_{0110})$$
$$=2(\log(z)+2p\pi i,\log(1-z)+2q\pi i)+2(\log(\oz)-2p\pi i,\log(1-\oz)-2q\pi i).$$

\end{cor}
\begin{pf}This follows from \hyperref[corsum]{Corollary \ref*{corsum}} and \hyperref[corfive]{Corollary \ref*{corfive}}, namely
the right hand side will be obtained by plugging
twice the right hand side of \hyperref[corfive]{Corollary \ref*{corfive}} in place of the second, third and fourth summand of the 
right hand side of \hyperref[corsum]{Corollary \ref*{corsum}}. The result then follows
by an elementary 
computation in each of the seven cases arising by the different case distinctions. Let us spell it out for the cases 
$0<Re(z)<1$ and $Im(z)>0$ or $Im(z)<0$, respectively.

If $0<Re(z)<1$ and $Im(z)>0$, then we obtain from \hyperref[corsum]{Corollary \ref*{corsum}} and \hyperref[corfive]{Corollary \ref*{corfive}} that the wanted
sum is $$2(\log(z)+2p\pi i,\log(1-z)+2q\pi i)+2(\log(\oz)-2p\pi i,\log(1-\oz)-2q\pi i)$$
$$+\chi(-\log(\oz)+\log(\frac{\oz-z}{\oz-1})+\log(1-\oz)-\log(1-\frac{\oz}{z})-\log(\frac{z}{\oz}))$$
and one easily checks that in the argument of $\chi$ everything cancels out such that the last summand is actually $\chi(0)=0$.

If $0<Re(z)<1$ and $Im(z)<0$, then the last summand will be
$$\chi(-\log(\oz)+2\log(\frac{(1-z)\oz}{\oz-z})+2(q-1)\pi i-2\log(\frac{1-z}{1-\oz})-\log(\frac{\oz-z}{\oz-1})-\log(1-\oz)+
\log(1-\frac{\oz}{z})+\log(\frac{z}{\oz})+2q\pi i)$$
and this time the cancellation occurs because of $$\chi(2\log(z-\oz)-2\log(\oz-z)-2\pi i)=\chi(-4\pi i)=0.$$

A very similar 
calculation works in the other five cases.

\end{pf}

\subsection{Proof of \hyperref[thm1]{Theorem \ref*{thm1}}}\label{proof}
{\em Let $M$ be a finite-volume, orientable, hyperbolic $3$-manifold.
Let $\tau\colon PSL(2,\C)\to SO(3,1)$ be the isomorphism $PSL(2,\C)\to SO(3,1)$.
Then for each boundary-unipotent representation $\rho\colon \pi_1M\to PSL(2,\C)$ one has 
$$(\tau\circ \rho)_*\left[M,\partial M\right]=2\rho_*\left[M,\partial M\right]+
2\overline{\rho_*\left[M,\partial M\right]}\in\widehat{\mathcal{B}}(\C)$$
if $\rho$ lifts to a boundary-unipotent representation $\pi_1M\to SL(2,\C)$ 
(in particular if $M$ is closed)
and $$(\tau \circ\rho)_*\left[M,\partial M\right]=2\rho_*\left[M,\partial M\right]+
2\overline{\rho_*\left[M,\partial M\right]}\in\widehat{\mathcal{B}}(\C)_{PSL}$$
otherwise.}\\

\begin{pf}
Let us first assume that $\rho\colon\pi_1M\to PSL(2,\C)$ lifts to a boundary-unipotent representation $\pi_1 M\to SL(2,\C)$, which abusing notation we will also denote by $\rho$. 

Fix some generalized ideal triangulation $M=\cup_{k=1}^r T_k$ whose lift to $\widetilde{M}$ admits a $\rho$-equivariant decoration whose ptolemy coordinates $c_t^k$ in the sense of \hyperref[ptolemy1]{Definition \ref*{ptolemy1}} are generic in the sense of \hyperref[generi]{Definition \ref*{generi}}. 
(Such a triangulation exists by \cite[Proposition 5.4]{gtz}.) For $k=1,\ldots,r$ let 
$\lambda(c^k)$ be defined as in \hyperref[invariant]{Definition \ref*{invariant}} (with $\alpha=0$) and define $z_k\in\C, p_q,q_k\in\Z$ via
$$\lambda(c^k)=(\log(z_k)+2p_k\pi i,\log(1-z_k)+2q_k\pi i).$$
By definition we have $$\rho_*\left[M,\partial M\right]=\sum_{k=1}^r\epsilon_k \lambda(c^k)$$
with the sign $\epsilon_k=\pm 1$ depending on orientation of $T_k$. By \hyperref[cancel]{Corollary \ref*{cancel}} we have 
$$(\rho\otimes\overline{\rho})_*\left[M,\partial M\right]=\sum_{k=1}^r\epsilon_k (2(\log(z_k)+2p_k\pi i,\log(1-z_k)+2q_k\pi i)+2(\log(\oz_k)-2p_k\pi i,\log(1-\oz_k)-2q_k\pi i)),$$
from which we conclude $$(\rho\otimes\overline{\rho})_*\left[M,\partial M\right]=2\rho_*\left[M,\partial M\right]+
2\overline{\rho_*\left[M,\partial M\right]}.$$

Finally it is known from the representation theory of $SL(2,\C)$ that the 2-fold covering $SL(2,\C)\to SO(3,1)$ is conjugate in $GL(4,\C)$ to $$id\otimes\overline{id}\colon SL(2,\C)\to SL(4,\C).$$ 
This implies that $\rho\otimes\overline{\rho}$ is conjugate to $\tau\circ \rho$, so we obtain the wanted equality
$$(\tau\circ\rho)_*\left[M,\partial M\right]=2\rho_*\left[M,\partial M\right]+
2\overline{\rho_*\left[M,\partial M\right]}.$$

If $\rho$ does not lift to $SL(2,\C)$, then the ptolemy coordinates $a,\ldots,f$ in \hyperref[standard]{Section \ref*{standard}} are only defined up to sign. We can choose some sign, so that for each simplex $T_k$ we have some equality of the kind $\pm a_kf_k\pm c_kd_k=b_ke_k$ (with certain signs) and can then still do all the computations in $\PC(\C)_{PSL}$ to get the equality there.\end{pf}\\

Remark: The proof of \hyperref[thm1]{Theorem \ref*{thm1}} via \hyperref[cancel]{Corollary \ref*{cancel}} might leave the impression that the equality from Theorem 1 already holds simplexwise, but one should be aware that this is just a (surprising) effect of the special choice of logarithms before \hyperref[ptolemy]{Lemma \ref*{ptolemy}}. With other (perhaps more natural) choices of logarithms the wanted equality would not hold simplexwise, rather additional contributions from different simplices would cancel out.

\subsection{Non-liftable PSL(2,C)-representations}\label{nonlift}

When a boundary-unipotent representation $\rho\colon\pi_1M\to PSL(2,\C)$ does not lift\footnotemark\footnotetext{This is in particular the case for the hyperbolic monodromy of a cusped hyperbolic $3$-manifold. Although such representations by \cite{cul} can always be lifted to $SL(2,\C)$, it is proved in \cite{cal} that boundaries of incompressible surfaces necessarily lift to parabolic elements with eigenvalue $-1$. Since these manifolds are always Haken this implies that they can not have a boundary-unipotent lift to $SL(2,\C)$. See \cite[Proposition 9.20]{gtz}.} to a boundary-unipotent representation $\pi_1M\to SL(2,\C)$, then its Chern-Simons invariant is only defined modulo $\pi^2$ and so of course the equality in $CS(\tau\circ\rho)=4CS(\rho)$ can only hold modulo $\pi^2$. 

However, since $\rho\otimes\overline{\rho}$ is well-defined as a boundary-unipotent representation to $SL(4,\C)$ it actually makes sense to compute its Chern-Simons invariant modulo $4\pi^2$ and we will describe in this section how to do this calculation.

Given a boundary-unipotent representation $\pi_1M\to PSL(2,\C)$ its obstruction to lifting it 
as a boundary-unipotent representation to $SL(2,\C)$ is represented by a $2$-cycle $\sigma\in Z_2(K,\partial K;\Z/2\Z)$. 
Depending on $\sigma$ the Ptolemy coordinates have to satisfy a certain simplexwise equation. Namely if for 
$i=0,\ldots,3$ we denote by $\sigma_i$ the value of $\sigma$ on the face opposite to the $i$-th vertex, then 
$$\sigma_0\sigma_1 af+\sigma_0\sigma_3 cd=\sigma_0\sigma_2 be,$$
 see \cite[Definition 3.5]{ggz}. (Now we fix one simplex and denote its ptolemy coordinates for the $SL(2,\C)$-representation again by $a,\ldots,f$ as in \hyperref[standard]{Section \ref*{standard}}.)

For a given $\sigma$ one can then use the methods from \cite[Section 9]{gtz} to simplexwise compute decorations for a simplex with given ptolemy coordinates $a=c_{01},b=c_{02},\ldots,f=c_{23}$.
(According to \cite[Definition 9.23]{gtz} the "diamond coordinates" of a face have to be multiplied by a sign according to its appearance in the obstruction cycle.) The result of the computations is in the end that a decoration of a $3$-simplex is given by

$$((\begin{array}{cc}1&0\\
0&1\end{array})N
,(
\begin{array}{cc}0&-\frac{1}{a}\\
a&0\end{array})N
,(\begin{array}{cc}\frac{- d}{a}&\mp\frac{1}{b}\\
\pm b&0\end{array})N
,(\begin{array}{cc}\frac{- e}{a}&\mp\frac{1}{c}\\
\pm c&0\end{array})N),
$$ 
where the sign in front of $b$ is positive if and only if $\sigma_3=0$, and the sign in front of $c$ is positive if and only if $\sigma_2=0$. 

When the signs in $\pm b$ and $\pm c$ are chosen as above, then we have an equality $a(\pm f)+ (\pm c)d=(\pm b)e$ with the sign in front of $f$ being positive if and only if $\sigma_1=0$.

This means that we can do the computations from \hyperref[tensor]{Section \ref*{tensor}} but with $b,c,f$ replaced by $\pm b,\pm c,\pm f$ according to the values of the obstruction cycle. (Note that the $f$-coordinate does not appear in the decoration, but it made its entrance in the calculations of \hyperref[tensor]{Section \ref*{tensor}} indirectly through the formula $af+cd=be$. For this reason we also have to change the sign of $f$ accordingly.)


Then one can use \hyperref[ptolemy]{Lemma \ref*{ptolemy}} to compute 
$(\rho\otimes\overline{\rho})_*\left[M,\partial M\right]$ as an element in $\widehat{\mathcal{B}}(\C)$ (rather just in $\widehat{\mathcal{B}}(\C)_{PSL}$) and the Chern-Simons invariant modulo $4\pi^2$ (rather just modulo $\pi^2$).

As an illustration let us take up the example of the figure eight knot complement $S^3\setminus K(\frac{2}{5})$. Our computation in  \hyperref[figureeight]{Section \ref*{figureeight}} actually illustrates the general principle from this subsection. The result was that $(\rho\otimes\overline{\rho})_*\left[M,\partial M\right]=0$ holds in $\PC(\C)$, and not just in $\PC(\C)_{PSL}$ as it would result from \hyperref[thm1]{Theorem \ref*{thm1}}. Hence 
$$CS(S^3\setminus K(\frac{2}{5}),\tau\circ\iota)=0$$
holds even modulo $4\pi^2$ and not just modulo $\pi^2$.

\section{On components of character varieties}\label{compo}
This section has two purposes. In \hyperref[4r]{Section \ref*{4r}} we 
use \hyperref[thm1]{Theorem \ref*{thm1}} to derive 
\hyperref[3compo]{Corollary \ref*{3compo}} (and 
we discuss  
why we need Chern-Simons invariants to distinguish components of the character variety while the known local rigidity results 
do not suffice). In \hyperref[galo]{Section \ref*{galo}} and \hyperref[epimo]{Section \ref*{epimo}} we 
discuss two methods to construct more components of the character variety for specific classes of $3$-manifolds. In particular in \hyperref[epimo]{Section \ref*{epimo}} we prove \hyperref[arbi]{Theorem \ref*{arbi}} which provides esamples of knot complements with an arbitrarily large number of connected components of the character variety that are not distinguishable by Chern-Simons invariants.

To repeat the basic notions, the variety of representations $Hom(\Gamma,G)$ of a finitely generated group $\Gamma$ into an algebraic group $G$ is by definition the variety defined by the relations between the given generators. The character variety is its quotient $Hom(\Gamma,G)//G$ in the sense of geometric invariant theory. We consider these varieties with the euclidean topology (not the Zariski topology). If $G$ is connected, then connected components of the character variety correspond to connected components of the representation variety. By a 
result of Goldman, the number of connected components is finite if $\Gamma$ is finitely generated and $G$ semisimple with finite center.

For closed manifolds we will consider the variety of characters of all representations $$X(\pi_1M,SL(n,\C))=Hom(\pi_1M,SL(n,\C))//SL(n,\C).$$ For manifolds with boundary we will consider only those representations whose restriction to $\partial M$ has unipotent image:
$$X_{bup}(\pi_1M,SL(n,\C))=\left\{\rho\in Hom(\pi_1M,SL(n,\C))\colon \rho(\pi_1\partial M)\subset N\right\}//SL(n,\C).$$
Volume and Chern-Simons invariant are constant on connected components of these character varieties, see \cite[Theorem 3.4]{lod} or \cite[Proposition 4.1]{bas}.

\subsection{SL(4,R)-representations factoring over PSL(2,C)}\label{4r}

Any representation of the form $\rho\otimes\overline{\rho}$ (for a representation $\rho$ into $SL(2,\C)$) can be conjugated into $SL(4,\R)$. 
Indeed it is easy to check that the matrix entries of $\rho\otimes\overline{\rho}$ with respect to the basis 

$\left\{e_1\otimes e_1,e_1\otimes e_2+e_2\otimes e_1,i(e_1\otimes e_2-e_2\otimes e_1),e_2\otimes e_2\right\}$ are all real. 

Even better, there is an isomorphism $PSL(2,\C)\to SO(3,1)$ explicitly defined by $\tau\left(\begin{array}{cc}a&b\\
c&d\end{array}\right)=1/2\ \times$
$$\hspace{-0.5in}\left(\begin{array}{cccc}\mid a\mid^2+\mid b\mid^2+\mid c\mid^2+\mid d\mid^2&\oa b+a\ob+\oc d+c\od&i(\oa b-a\ob+\oc d-c\od)& \mid a\mid^2-\mid b\mid^2+\mid c\mid^2-\mid d\mid^2\\
a\oc+\oa c+b\od +\ob d&a\od+\oa d+b\oc+\ob c&i(\oa d+b\oc -a\od -\ob c)&a\oc+\oa c-b\od -\ob d\\
-i(\ob d+\oa c-a \oc-b\od)&-i(a\od +\ob c-\oa d-b\oc)&a\od+\oa d-b\oc-\ob c&-i(\oa c+b\od -a\oc -\ob d)\\
\mid a\mid^2+\mid b\mid^2-\mid c\mid^2-\mid d\mid^2&a\ob+\oa b-c\od-\oc d&i(\oa b+c\od -a\ob -\oc d)& \mid a\mid^2-\mid b\mid^2-\mid c\mid^2+\mid d\mid^2\end{array}\right).$$
and it is known from the representation theory of the Lorentz group that the "four-vector representation", i.e., the corresponding 2-fold covering $SL(2,\C)\to SO(3,1)$, is equivalent to the representation $\rho_{1,1}=id\otimes\overline{id}$. So for any representation $\rho$ we have $$\tau\circ\rho\sim\rho\otimes\overline{\rho}.$$
(That was why our computations in \hyperref[compu]{Section \ref*{compu}} also implied $Vol(\tau\circ\rho)=0, CS(\tau\circ\rho)=4CS(\rho)$.) 

We remark that this isomorphism is not well-behaved with respect to genericity of ptolemy coordinates in the sense of \hyperref[ptolemy1]{Definition \ref*{ptolemy1}}. Given a triangulation and a generic $PSL(2,\C)$-decoration one may well get a non-generic $SO(3,1)$-decoration after applying the isomorphism. For this reason it was more convenient to compute Chern-Simons invariants for $\rho\otimes\overline{\rho}$ as we did in \hyperref[compu]{Section \ref*{compu}}, rather than trying to compute them for $\tau\circ\rho$ directly.

Besides the trivial representation and the four-vector representation there is only one more representation of $SL(2,\C)$ in $SL(4,\R)$ namely
the representation $\kappa:SL(2,\C)\rightarrow SL(4,\R)\subset SL(4,\C)$ defined by
$$\kappa\left(\begin{array}{cc}a_1+a_2i&b_1+b_2i\\
c_1+c_2i&d_1+d_2i\end{array}\right)=\left(\begin{array}{cccc}a_1&a_2&b_1&b_2\\
-a_2&a_1&-b_2&b_1\\
c_1&c_2&d_1&d_2\\
-c_2&c_1&-d_2&d_1\end{array}\right).$$
Again, application of this representation to the standard decoration would yield a nongeneric decoration. However conjugation with the matrix $$\left(\begin{array}{cccc}1&0&1&0\\
i&0&-i&0\\
0&1&0&1\\
0&i&0&-i\end{array}\right)\in GL(4,\C)$$ shows that $\kappa$ is equivalent to
$id\oplus\overline{id}$, which then with \hyperref[additivity]{Lemma \ref*{additivity}}
implies that $$Vol(\kappa\circ\iota)=0,\hspace{0.2in}CS(\kappa\circ\iota)=2CS(M).$$

Let $M$ be a finite-volume hyperbolic $3$-manifold and $\iota\colon\pi_1M \to SL(2,\C)$
a lift of the hyperbolic monodromy $\iota\colon\pi_1M \to SL(2,\C)$. To the best of my knowledge it has not been known so far whether the three characters of $\tau\circ\iota$, $\kappa\circ\iota$ and the trivial representation necessarily belong to different components in the $SL(4,\R)$-character variety. It is known that in general neither of this representations has to be locally rigid (i.e., the corresponding component of the character variety is not just an isolated point).
For example, \cite{klw} constructs (for some $3$-manifolds) a 1-parameter deformation of $\kappa\circ \iota$ in $Sp(4,\R)\subset SL(4,\R)$, which is therefore not an isolated point in the character variety. It is also known that the representations $\iota\otimes \overline{\iota}=(\rho_2\otimes\overline{\rho}_2)\circ\iota$ are not always locally rigid, for example \cite{clt} shows that local rigidity in $SL(4,\R)$ does not hold for exactly 52 of the first 4500 closed, orientable, hyperbolic $3$-manifolds with 2-generator fundamental group in the Hodgson-Weeks census.
And for the trivial representation $\nu\colon\Gamma\to GL(n,\C)$ 
one has $H^1(\Gamma,Ad(\nu))=H^1(\Gamma,\C^n),$ i.e., one has deformations corresponding to homomorphisms $\pi_1M\to H_1M\to\C^n$.


But for hyperbolic $3$-manifolds with $CS(M)\not=0$ we can now use $$CS(\tau\circ\iota)=4CS(M)\not=0\not=2CS(M)=CS(\kappa\circ\iota)$$ to obtain that $\tau\circ\iota$, $\kappa\circ\iota$ and the trivial character belong to three different components of the character variety. Hence for hyperbolic $3$-manifolds with $CS(M)\not=0$ the $SL(4,\R)$-character variety has at least three connected components, which proves \hyperref[3compo]{Corollary \ref*{3compo}}.

Similarly we obtain in this case at least 10 components in the $SL(4,\C)$-character variety, as we had explained in the introduction.

For the figure eight knot complement, the experimental results from \cite{pto} suggest that there are not more than three components of the $SL(4,\R)$-character variety containing characters of irreducible representations. 
In particular, even though $CS(M)=0$ in this case, the experimental result seems to suggest that $\tau\circ\iota$ and $\kappa\circ\iota$ are isolated points and so again we have 3 components.

For other knots one may have more than three components and
the following two subsections will shortly discuss two approaches to construct some of them.
\subsection{Using Galois actions}\label{galo}
One can frequently get more components by applying Galois actions. Namely, for each finite-volume hyperbolic $3$-manifold $M$, the image of the hyperbolic monodromy $\rho\colon \pi_1M\to PSL(2,\C)$ is contained in $PSL(2,K)$ for some number field $K$ and then any element of the Galois group $\sigma\in Gal(K:\Q)$ provides a (non-discrete) representation $$\rho^\sigma\colon\pi_1M\to PSL(2,K)\subset PSL(2,\C).$$ Composition with representations $PSL(2,\C)\to SL(m,\C)$ then produces more $SL(m,\C)$-representations of $\pi_1M$. One should also note that local rigidity results from \cite{mp} and \cite{rig} carry over to the Galois conjugate representations because the group cohomology $H^1(\Gamma,Ad(\rho))$ is compatible with Galois conjugations.

As an illustration let us look at the 2-bridge knot $K(\frac{7}{3})$, which is the knot considered in \cite[Example 10.1]{gtz}. Its hyperbolic monodromy has image in $PSL(2,\Q(x))$ with $x$ the unique root of positive imaginary part of $$x^3-x^2+1=0.$$ 
Letting $x^\prime$ be the real root, the Galois automorphism $\Q(x)\to \Q(x^\prime)$ yields a representation $\rho\colon\Gamma \to PSL(2,\R)$, which according to \cite[Example 6.16]{zic1} has Chern-Simons invariant $-1.1134\ldots$, while the hyperbolic monodromy and its complex conjugate have Chern-Simons invariants $\pm 3.0241\ldots$. Then $\rho\otimes\overline{\rho}$ is a representation to $SL(4,\R)$ of Chern-Simons invariant $-4.453\ldots$ which consequently does not belong to any of the three other components.\footnotemark\footnotetext{I do not know whether it has some meaning that this value coincides with the Chern-Simons invariants of one of the $SL(3,\R)$-representations in \cite[Example 10.1]{gtz}.}

\subsection{Using epimorphisms of knot groups}\label{epimo}
Recall that for each rational number $\frac{p}{q}$ one has a 2-bridge link $K(\frac{p}{q})$ which is determined by the coefficients $b_1,b_2,\ldots,b_k$ in the continued fraction expansion $$\frac{p}{q}=\left[b_1,b_2,\ldots,b_k\right]=\frac{1}{b_1+\frac{1}{b_2+\frac{1}{b_3+\ldots}}}$$ and has a plane diagram which (for even $k$) looks as in he picture below, where the number $\pm b_i\in\Z$ inscribed in a box means the number of half-twists in that box with positive numbers corresponding to right-handed half-twists and negative numbers corresponding to left-handed ones. 

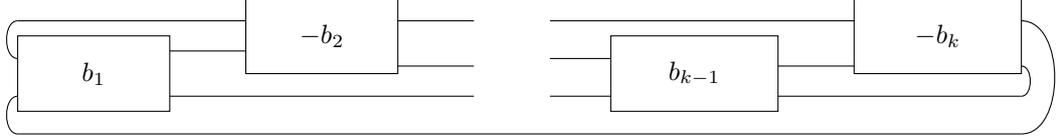
\begin{tikzpicture}

\draw(0,0.5)rectangle(2,1.5);
\node at (1,1){$b_1$};
\draw(2,0.7)--(6,0.7);
\draw(2,1.3)--(3,1.3);

\draw(3,1)rectangle(5,2);\node at (4,1.5){$-b_2$};
\draw(5,1.7)--(6,1.7);
\draw(5,1.1)--(6,1.1);
\draw(7,1.2)--(7.8,1.2);\draw(7,0.7)--(7.8,0.7);
\draw(7.8,0.5)rectangle(10,1.5);\node at (8.9,1){$b_{k-1}$};
\draw(7,1.7)--(11,1.7);
\draw(10,1.1)--(11,1.1);
\draw(11,1)rectangle(13.2,2);\node at (12.1,1.5){$-b_k$};
\draw(10,0.7)--(13.2,0.7);
\draw(13.2,1.1) to[out=0,in=360] (13.2,0.7);

\draw(0,0.2)--(13.2,0.2);
\draw(13.2,1.7) to[out=0,in=360] (13.2,0.2);
\draw(0,0.7) to [out=180,in=180] (0,0.2);
\draw(0,1.7)--(3,1.7);
\draw(0,1.7) to [out=180,in=180] (0,1.2);

\end{tikzpicture}

The link $K(\frac{p}{q})$ is hyperbolic for $p>1$.

We will use the following result from \cite{ors} which is a combination of Proposition 5.1 (or the equivalent Proposition 5.2) and Theorem 6.1 in that paper. Following \cite[Section 5]{ors} we will use the notation ${\mathbf b}=(b_1,\ldots,b_k), {\mathbf b}^{-1}=(b_k,\ldots,b_1), \epsilon{\mathbf b}=(\epsilon b_1,\ldots,\epsilon b_k), \epsilon{\mathbf b}^{-1}=(\epsilon b_k,\ldots,\epsilon b_1)$ with $\epsilon\in\left\{1,-1\right\}$.
\begin{pro}\label{orsresult}(Ohtsuki-Riley-Sakuma): Assume $\frac{p}{q}=\left[b_1,b_2,\ldots,b_k\right]$ and 
$$\frac{\tilde{p}}{\tilde{q}}=2c+
\left[\epsilon_1{\mathbf b},2c_1,\epsilon_2{\mathbf b}^{-1},2c_2,\epsilon_3{\mathbf b},\ldots,2c_{l-1},\epsilon_l{\mathbf b}^{(-1)^{l-1}}\right]$$ 
for\footnotemark\footnotetext{The statement in \cite[Proposition 5.1]{ors} assumes $c$ positive, but the equivalent \cite[Proposition 5.2]{ors} and its proof do not require this assumption.} $l\in\N,c\in\N\cup\left\{0\right\}$ and some $\epsilon_1,\ldots,\epsilon_l\in\left\{1,-1\right\},c_1,\ldots,c_{l-1}\in\Z$. 
Then there exists a proper, branched fold map $f\colon (S^3,K(\frac{p}{q}))\to (S^3,K(\frac{\tilde{p}}{\tilde{q}}))$ inducing an epimorphism $$f_*\colon\pi_1(S^3\setminus K(\frac{p}{q}))\to \pi_1(S^3\setminus K(\frac{\tilde{p}}{\tilde{q}})).$$ \end{pro}

The map sends a meridian to a meridian and a longitude to the power of a longitude (see \cite[Section 4]{hs}), in particular it is boundary-preserving and $f^*$ defines a map between boundary-unipotent character varieties.

Using this result of theirs, Ohtsuki, Riley and Sakuma proved in \cite[Corollary 7.3]{ors} that there are 2-bridge link character varieties with arbitrarily large number of irreducible components.\footnotemark\footnotetext{They stated this theorem for $SL(2,\C)$, but their argument also applies to $SL(m,\C)$-character varieties for $m>2$. One just has to consider the geometric representation (that is, the composition of the lifted hyperbolic monodromy $\rho_0$ with the irreducible representation $SL(2,\C)\to SL(m,\C)$) in place of $\rho_0$ to adapt their proof.} 
Namely they constructed a sequence of hyperbolic links $K_n$ with non-injective epimorphisms $$f_n\colon \pi_1(S^3\setminus K_n)\to \pi_1(S^3\setminus K_{n-1}).$$ Let $C_n$ be the 
irreducible component containing the defining representation for the hyperbolic metric of $S^3\setminus K_n$, because the latter is an irreducible representation. The epimorphism $f_n$ induces an injective, regular map $f_n^*$ of the representation varieties. Then $f_n^*C_{n-1}$ is an irreducible component consisting of non-faithful 
representations and thus does not agree with the irreducible component $C_n$ containing the defining representation for the hyperbolic metric of $S^3\setminus K_{n-1}$. Iterating this 
argument they get at least $n$ irreducible components for the character variety of $K_n$: the argument above shows that the irreducible components $f_i^*C_{i-1}$ and $C_i$ are distinct, 
and this implies (since preimages of irreducible sets under injective, regular maps are irreducible)
that 
the irreducible components $f_n^*\ldots f_i^*C_{i-1}$ and $f_n^*\ldots f_{i+1}^*C_i$ are distinct.

Since we are considering the euclidean topology rather than the Zariski topology, a connected component may consist of several irreducible components. So having many irreducible components does, a priori, say nothing about the number of components of the character variety. The following argument shall show, however, that the Ohtsuki-Riley-Sakuma construction actually can yield 2-bridge link character varieties with an arbitrarily large number of connected components.

\begin{pro}\label{2br}There are $2$-bridge knots with the boundary-unipotent $SL(m,\C)$-character variety $X_{bup}$ having arbitrarily large number of connected components.\end{pro}
\begin{pf}

As a special instance of the construction in \cite{ors} let us consider the 
sequence $L_n=K(\frac{p_n}{q_n})$ of 2-bridge knots such that $\frac{p_n}{q_n}$ has a continued fraction expansion 
$$\frac{p_n}{q_n}=\left[2,2,2,2,\ldots,2,2\right]$$ with $2\cdot3^{n-1}$ coefficients all equal to $2$.

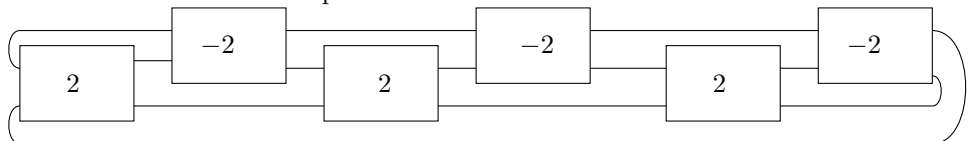
\begin{tikzpicture}

\draw(0,0.5)rectangle(1.5,1.5);
\node at (0.7,1){$2$};
\draw(1.5,0.7)--(4,0.7);
\draw(1.5,1.3)--(2,1.3);

\draw(2,1)rectangle(3.5,2);
\node at (2.6,1.5){$-2$};
\draw(3.5,1.7)--(6,1.7);
\draw(3.5,1.2)--(4,1.2);
\draw(4,0.5)rectangle(5.5,1.5);
\node at (4.8,1){$2$};
\draw(6,1)rectangle(7.5,2);
\node at (6.8,1.5){$-2$};
\draw(5.5,1.2)--(6,1.2);
\draw(7.5,1.2)--(8.5,1.2);

\draw(5.5,0.7)--(8.5,0.7);
\draw(8.5,0.5)rectangle(10,1.5);
\node at (9.2,1){$2$};
\draw(7.5,1.7)--(10.5,1.7);
\draw(10,1.2)--(10.5,1.2);
\draw(10.5,1)rectangle(12,2);
\node at (11.1,1.5){$-2$};
\draw(10,0.7)--(12,0.7);
\draw(12,1.1) to[out=0,in=360] (12,0.7);

\draw(0,0.2)--(12,0.2);
\draw(12,1.7) to[out=0,in=360] (12,0.2);

\draw(0,0.7) to [out=180,in=180] (0,0.2);
\draw(0,1.7)--(2,1.7);

\draw(0,1.7) to [out=180,in=180] (0,1.2);

\end{tikzpicture}

So $L_1=K(\frac{2}{5})$ is the figure eight knot, $L_2$ (shown in the picture above) is $K(\frac{70}{169})$
and so on. Then, for any $n\ge 2$, \hyperref[orsresult]{Proposition \ref*{orsresult}} applied with $c=0, l=3, c_1=c_2=0, \epsilon_1=1=\epsilon_2=\epsilon_3=1$ and ${\mathbf b}={\mathbf b}^{-1}=\left[2.\ldots,2\right]$ provides us with boundary-preserving maps $$f_n\colon S^3\setminus L_n\to S^3\setminus L_{n-1}.$$ 

By \cite[Proposition 6.2]{ors} the degree of the map $f$ constructed in \hyperref[orsresult]{Proposition \ref*{orsresult}} can be computed as $deg(f)=\sum_{j=1}^l(-1)^{j+1}\epsilon_j$. So in our example we have $deg(f_n)=1$.


We now want to prove that pulling back components via $f_n^*$ one obtains many components that can be distinguished by volume.

There are linear bounds $$C_1 tw(L)-C_2\le vol(S^3\setminus L)\le C_3 tw(L)+C_4$$ (with explicit constants $C_1,C_2,C_3,C_4$) for the hyperbolic volume in terms of the twist number of a prime alternating diagram, see \cite{lac},\cite{fg}. We can apply this to our diagrams because they are alternating. (This is a consequence of all the coefficients $b_i$ having the same sign and can of course also be seen directly from the link diagram.)
The number of twist regions for the canonical diagram of $L_n$ is 
$$tw(L_n)=2\cdot3^{n-1}.$$  
In particular, in our situation $vol(S^3-L_n)$ will go to infinity exponentially. 

Let $\rho_n\colon\pi_1(S^3\setminus L_n)\to SL(m,\C)$ be 
the geometric representation. By \hyperref[invpol]{Lemma \ref*{invpol}} 
we have $$vol(\rho_n)=\frac{(m+1)m(m-1)}{6}Vol(S^3\setminus L_n),$$
so also $vol(\rho_n)$ grows exponentially in $n$. In particular, 
there are natural numbers $N_1,N_2$ such that for all $N\ge N_1$ there are at least $N-N_2$ distinct values among the volumes $$vol(\rho_n), n=1,\ldots,N.$$

Let $n\le N$. Denote by $f_{k*}$ the homomorphism of fundamental groups induced from $f_k$. From $deg(f_N)=\ldots=deg(f_{n+1})=1$ we obtain $$vol(\rho_n f_{n+1,*}\ldots f_{N,*})=vol(\rho_n),$$
see \hyperref[degree]{Lemma \ref*{degree}}. Hence there are also at least $N-N_2$ distinct values among the volumes 
$$ vol(\rho_n f_{n+1,*}\ldots f_{N,*}), n=1,\ldots,N$$
for $N\ge N_1$. In particular the $SL(m,\C)$-character variety of $S^3\setminus L_N$ has at least $N-N_2$ connected components.
\end{pf}\\

With the help of the following proposition we can even show that the above constructed components all correspond to representations of trivial Chern-Simons invariant. 
\begin{pro}\label{2brcs}\footnotemark\footnotetext{Experimental evidence from the knot table \cite{cl} suggests that also the converse of this proposition might be true.}The complement of a 2-bridge knot $K(\frac{p}{q})$ has vanishing $PSL(2,\C)$-Chern-Simons invariant if the continued fraction expansion $\frac{p}{q}=\left[a_1,\ldots,a_k\right]$ is symmetric in the sense that $$a_1=a_k,a_2=a_{k-1},a_3=a_{k-2},\ldots$$\end{pro}
\begin{pf}(\cite{mathov}) The symmetry of the continued fraction expansion is equivalent to $q^2\equiv -1\ mod\ p$, see \cite{ser}, or \cite{smi} for a more general statement. If $K(\frac{p}{q})$ is a knot, then $q$ is odd and we obtain $q^2\equiv -1\ mod\ 2p$. From the "Korollar zu Satz 4" in \cite{sch} this is equivalent to $K(\frac{p}{q})$ being amphicheiral. According to \cite{mo} this implies $CS(S^3\setminus K(\frac{p}{q}))=0$.
\end{pf}\\

This shows that all the examples in the proof of \hyperref[2br]{Proposition \ref*{2br}} will have vanishing Chern-Simons invariant and so \hyperref[2br]{Proposition \ref*{2br}} actually produces examples of manifolds with $SL(m,\C)$-character variety having arbitrarily large number of connected components of vanishing Chern-Simons invariants. This proves \hyperref[arbi]{Theorem \ref*{arbi}} from the introduction.

The argument for \hyperref[2br]{Proposition \ref*{2br}} does not apply to $SL(m,\R)$-character varieties because the volume of $SL(m,\R)$-representations is always zero. Instead one should use Chern-Simons invariants, but of course the above argument using growth of volumes does not adapt because no such statement can be true for Chern-Simons invariants. Still one can use explicitly computed values to construct distinct components in specific examples. E.g.\ for the 2-bridge knot $K(\frac{3}{7})$ one can use the degree 2-map $$S^3\setminus K(\frac{23}{53})\to S^3\setminus K(\frac{3}{7})$$ (that one obtains from \hyperref[orsresult]{Proposition \ref*{orsresult}} with ${\mathbf b}=(2,3), l=2, c=c_1=0,\epsilon_1=1,\epsilon_2=-1$) and the inequality $$CS( S^3\setminus K(\frac{23}{53}))\not=2CS( S^3\setminus K(\frac{3}{7}))$$ to obtain an additional component in the $SL(4,\R)$-character variety of the 2-bridge knot $K(\frac{23}{53})$. A similar argument applies to other components of the character variety and to many other knots.

\end{document}